\documentclass[12pt]{amsart}
\usepackage[square,compress,comma, numbers,sort]{natbib}
\usepackage[colorlinks=true, citecolor=blue, linkcolor=blue]{hyperref}
\usepackage{amsfonts}
\allowdisplaybreaks[4]
\usepackage{graphicx}
\usepackage{epstopdf}
\usepackage{subfigure}
\usepackage{caption}
\usepackage[svgnames]{xcolor}
\usepackage{listings}
\usepackage{amssymb}
\usepackage{amsmath}

\usepackage{color}
\usepackage{enumerate}
\usepackage{float}
\newcommand{\abs}[1]{\left\lvert #1 \right\rvert}

\def\E#1{\mathbb{E}\left \{#1 \right\}}
\def\I#1{\mathbb{I}\left (#1 \right)}

\definecolor{c20}{rgb}{0.,0.7,0.}
\definecolor{c30}{rgb}{0.,0.,1.}
\definecolor{c40}{rgb}{1,0.1,0.7}
\definecolor{c50}{rgb}{1,0,0}
\definecolor{c60}{rgb}{1,0.9,0.1}
\definecolor{c70}{rgb}{0.50,1.00,0.00}
\def\cL#1{\textcolor{c30}{#1}}
\def\cL#1{#1}

\def\N{\mathbb{N}}
\numberwithin{equation}{section}
\newtheorem{theo}{Theorem}[section]
\newtheorem{prop}{Proposition}[section]
\newtheorem{sat}[theo]{Proposition}
\newtheorem{de}[theo]{Definition}
\newtheorem{lem}{Lemma}[section]

\newtheorem{korr}[theo]{Corollary}
\newtheorem{remark}[theo]{Remark}
\newtheorem{remarks}[theo]{Remarks}
\numberwithin{equation}{section}

\newcommand{\prooftheo}[1]{ \textsc{Proof of Theorem} \ref{#1} }
\newcommand{\proofprop}[1]{\textsc{Proof of Proposition} \ref{#1}}

\newcommand{\pk}[1]{\mathbb{P} \left\{ #1 \right\} }

\newcommand{\QED}{\hfill $\Box$}
\newcommand{\COM}[1]{}
\def\IF{\infty}
\newcommand{\R}{\mathbb{R}}
\newcommand{\inr}{\in \R}

\def\polhk#1{\setbox0=\hbox{#1}{\ooalign{\hidewidth
\lower1.5ex\hbox{`}\hidewidth\crcr\unhbox0}}}

\def\cL#1{\textcolor{c50}{#1}}
\def\cL#1{#1}

\def\zs#1{\textcolor{magenta}{#1}}
\def\zs#1{#1}

\def\rr#1{\textcolor{c50}{#1}}
\def\rr#1{#1}
\def\ree#1{\textcolor{c50}{#1}}
\def\ree#1{#1}

\newcommand{\sign}{\mathrm{sign}}

\usepackage{mathtools}

\newcommand{\neprop}[1]{{Proposition \ref{#1}}}
\newcommand{\netheo}[1]{{Theorem \ref{#1}}}

\def\IF{\infty}

\newcommand{\expon}[1]{\exp\left(#1\right)}

\def\Cov{\mathrm{Cov}}

\date{}

\newcommand{\BS}{\begin{sat}}
\newcommand{\ES}{\end{sat}}
\newcommand{\BT}{\begin{theo}}
\newcommand{\ET}{\end{theo}}
\newcommand{\BP}{\begin{prop}}
\newcommand{\EP}{\end{prop}}
\newcommand{\BK}{\begin{korr}}
\newcommand{\EK}{\end{korr}}

\newcommand{\BD}{\begin{de}}
\newcommand{\ED}{\end{de}}
\newcommand{\BIT}{\begin{itemize}}
\newcommand{\EIT}{\end{itemize}}
\newcommand{\BDI}{\begin{description}}
\newcommand{\EDI}{\end{description}}
\newcommand{\BRM}{\begin{remarks}}
\newcommand{\ERM}{\end{remarks}}
\newcommand{\BEL}{\begin{lem}}
\newcommand{\EEL}{\end{lem}}

\def\Cov{\mathrm{Cov}}

\def\MB{\mathcal{P\!B}}

\def\cL#1{\textcolor{c40}{#1}}
\def\cL#1{#1}
\def\CL#1{\textcolor{c40}{#1}}
\def\CL#1{#1}

  \def\td{\text{\rm d}}

\def\W{\mathcal{W}}
\def\N{\mathcal{N}}

\begin{document}
\title{On generalized Piterbarg-Berman function}

{\author{Chengxiu Ling}
     \address{Chengxiu Ling, Department of Mathematical Sciences, Xi'an Jiaotong-Liverpool University, Suzhou 215123, China}
     \email{chengxiu.ling@xjtlu.edu.cn}
}	

	{\author{Hong Zhang}
	\address{Hong Zhang, School of Mathematics and Statistics, Southwest University, Chongqing 400715, China}
	\email{crystal1994@email.swu.edu.cn}
}

{\author{Long Bai}
     \address{Long Bai, Department of Mathematical Sciences, Xi'an Jiaotong-Liverpool University, Suzhou 215123, China}
     \email{long.bai@xjtlu.edu.cn}
}	
	
{
	\bigskip
	
	\date{\today}
	\maketitle

\begin{abstract}
	\footnotesize{
	This paper aims to evaluate the Piterbarg-Berman function given by
	$$\MB_\alpha^h(x, E) = \int_\R e^z \pk{\int_E \mathbb{I}\left(\sqrt2B_\alpha(t) - |t|^\alpha - h(t) - z>0 \right) \td t > x} \td z,\quad x\in[0, {mes}(E)],$$
with  $h$ a drift function and $B_\alpha$ a fractional Brownian motion (fBm) with Hurst index $\alpha/2\in(0,1]$, i.e., a mean zero Gaussian process with continuous  sample paths and covariance function
\begin{align*}
\Cov(B_\alpha(s), B_\alpha(t)) = \frac12 (|s|^\alpha + |t|^\alpha - |s-t|^\alpha).
\end{align*}
This note specifies its explicit expression  for the fBms with $\alpha=1$ and $2$ when the drift function $h(t)=ct^\alpha, c>0$ and $E=\R_+\cup\{0\}$. For the Gaussian distribution $B_2$, we investigate $\MB_2^h(x, E)$ with general drift function $h(t)$ such that $h(t)+t^2$ being convex or concave, and finite interval $E=[a,b]$. Typical examples of $\MB_2^h(x, E)$ with $h(t)=c\abs{t}^\lambda-t^2$ and several bounds of $\MB_\alpha^h(x, E)$ are discussed. Numerical studies are carried out to illustrate all the findings.
	
 }


\end{abstract}

\section{Introduction}\label{Introduction}
Consider a centered Gaussian process $\{X(t),\, t \in \R\}$ with c\`{a}dl\`{a}d sample paths and let for $T>0$
\begin{align*}
 L_{u,T} := \int_{0}^T \mathbb{I}\left(X\left(t\right) >u\right) \td t
\end{align*}
be the sojourn time of $X$ above the level $u\inr$ during the observed period $[0,T]$, where $\mathbb I\left(\cdot\right)$ stands for the indicator function. In a series of papers culminating in \cite{Berman92Sojourns}, S. Berman derived results on the tail asymptotic \rr{behavior} of $\nu(u) L_{u, T}$ with  an appropriate scaling function $\nu(u)$  such that
\begin{align}\label{1.1}
p_{u, T}(x) = \pk{\nu(u) L_{u,T} > x} \sim C(x) \pk{\sup_{t\in [0, T]} X(t) > u}
\end{align}
as $u\to\IF$. This essentially builds a bridge of \cL{asymptotic}
behavior of the sojourn time $L_{u, T}$ and the extremal analysis of the Gaussian processes via the link function $C(x)$. However, the asymptotic function $C$ is in general difficult to \rr{obtain}  except very few special processes and approximations have been suggested to
evaluate. A related work is given by \cite{Akahori1995Some} for a standard Brownian motion with linear drift function. For a stationary and standard Gaussian process $X$ with correlation function $\rho$ satisfies the Pickands' assumption $\rho(t) = 1-|t|^\alpha[1+o(1)], \alpha\in(0, 2]$ for small $|t|$, \cite{Berman92Sojourns} showed an explicit form of function $C$ via the following tail distribution (see Theorem 3.3.1 therein)
\begin{align}\label{Alike_G}
G(x) = \pk{\int_\R \mathbb{I}(\sqrt2B_\alpha(t) - |t|^\alpha + \W >0) \td t > x},\quad x\ge0,
\end{align}
where $\W$ is a standard exponential distributed random variable, \rr{independent} of $B_\alpha $, a fractional Brownian motion  (fBm) with Hurst index $\alpha/2 \in(0,1]$, i.e., a mean zero Gaussian process with continuous  sample paths and covariance function
\begin{align*}
\Cov(B_\alpha(s), B_\alpha(t)) = \frac12 (|\ree{s}|^\alpha + |t|^\alpha - |s-t|^\alpha).
\end{align*}
The recent contribution \cite{debicki2017approximation} discussed   \eqref{1.1}  and gave the approximations of the related sojourn time of discrete form for locally stationary Gaussian processes, and \cite{debicki2018sojourn} investigated general Gaussian processes with strictly positive drift function. For more related discussions on ruin time and the extremal analysis of Gaussian processes and random fields in financial and insurance framework, we refer to \cite{Parisian2017, ParisianInfinite2018, Li2013The, Loeffen2014Occupation,D2002Ruin, H2004On,Dieker2005Extremes, D2017Uniform, debicki2015comparison, ling2018extremes}.

{Motivated} by the importance of the crucial function arising in the extremal behavior of  the sojourn time and the random processes involved, this paper studies thus a general form of function $G$ given in \eqref{Alike_G}. {Precisely},
define for a compact set $E$ in $\R$ and a continuous drift function $h$ on $E$
\begin{align}\label{MB_E}
\qquad \MB_{\alpha}^h(x,E):= \int_{\R} e^{z} \pk{ \int_E \mathbb I(\sqrt2 B_\alpha(t)  - |t|^\alpha - h(t) - z >0) \td t >x }  \td z,\quad x\in[0, \mathrm{mes}(E)]
\end{align}
and
\begin{align}\label{B_hat}
\MB_{\alpha}^h(x):= \lim_{T\to\IF}\frac{\MB_{\alpha}^h(x,[0,T])}{T^{\I{h=0}}},
\quad x \ge0
\end{align}
provided that the above integral and limits exist. For $h=0$, we suppress the superscript and write $\MB_\alpha(x)$ or $\MB_\alpha(x, E)$. Typical examples of the function $\MB_{\alpha}^h(x,E)$ {can be found} in \cite{debicki2017approximation} for $h=0, E=\R$, and \cite{debicki2018sojourn} for polynomial function $h$ and general \cL{interval} $E$.

\cL{Clearly, our setting is {very} common since  $\MB_{\alpha}(0)$ is simply the Pickands' constant, which values are known only for $\alpha =1, 2$, i.e., $\MB_1(0) = 1, \MB_2(0) =1/\sqrt\pi$, see e.g., \cite{piterbarg2012asymptotic, debicki2017generalized, harper2017pickands,PicandsA,Pit72, debicki2002ruin,DE2014,DiekerY,DEJ14,Pit20, Tabis, DM} for related studies on its expressions and bounds, while $\MB_{\alpha}^h(0)$ reduces to the  Piterbarg constants for strictly positive drift function.} The recent contribution  \cite{Long2018On} studied the basic properties of the generalized Piterbarg constant ${\MB}_{\alpha}^h(0)$  for power drift function, which are available for all $\alpha \in(0, 2]$. For more general \cL{studies} on  sojourn sets with moving boundary of the processes involved and applications in physics and finance fields, we refer to \cite{seuret2019sojourn} and among others.


The first result below is concerned with the explicit expression of $\MB_\alpha^h(x, E)$ for the standard Brownian motion and Gaussian distribution, i.e., the fBm with $\alpha=1$ and {$\alpha=2$}. Here, we focus simply on $E=\R_+ \cup\{0\}$ and positive drift function $h(t) = ct^\alpha, t\ge0, c>0$. In what follows, let  $\Psi(\cdot)$ and $\varphi(\cdot)$  be the survival function and probability density function of $\N\sim N(0,1)$, respectively.

\BT\label{T_one_two}
Let ${\MB}_\alpha^h(x)$ be the Piterbarg-Berman function given by \eqref{MB_E} with drift function $h(t) = ct^\alpha, t\ge0, c>0$. We have with $x_{c} = (1+c) \sqrt{x/2}$ and $x_{c}' = (1-c)\sqrt{x/2}$
\begin{align*}
{\MB}_{1}^h(x)  	
	&=\frac{(1+c)^2}c\Psi(x_c) - (1-c)x_c\varphi(x_c)  +\left[ \frac{(1-c)^2}{4(1+c)}x_c^2 +\frac{1+c}2\right]e^{-x_{c}^2/2}\\
	&\quad -(1+c) e^{-\cL{c}x} \left[\frac{1-c}c \Psi(x_{c}')- x_c'\varphi(x_c') + \frac{1+x_c'^2/2} 2 e^{-x_c'^2/2}\right]
\end{align*}
and
\begin{align*}
{\MB}_{2}^h(x)  	
	= \sqrt{\frac{1+c}{c}}\Psi\left( \sqrt {\frac{c(1+c)}2}x \right)e^{-\frac{(1+c)x^2}{4}} + \Psi\left( \frac{1-c}{\sqrt{2}}x \right)e^{-c x^2} - \Psi\left( \frac{1+c}{\sqrt{2}} x\right), \quad x\ge0.
\end{align*}
\ET
{\remark The explicit expression of ${\MB}_{1}^h(x)$ is obtained by the considerable analysis of  the stopping time and the random sojourn time involved:
\begin{align}\label{def_sojourn_BM}
\tau_z = \inf\{t\ge0: \sqrt2 B_1(t)- t - h(t) \ge z\},\quad Y_z = \int_0^\IF \mathbb I(\sqrt2 B_1(t) - t - h(t) -z>0) \td t,
\end{align}
which nice properties are referred to \cite{Borodin2012} due to the linear drift function. The  general case with non-linear drift function and finite time interval is an open question and it may require  definite efforts to develop the initial properties of $\tau_z$ and $Y_z$ involved.
}

The main methodology for the establishment of ${\MB}_{2}^h(x)$ above is essentially determined by the convex curve family $f_\N(t) = f_\N(t,z), z\inr$ (recall \eqref{MB_E})
\begin{align}\label{def_f}
f_\N(t) = h(t) + t^2 - \sqrt 2 \N t + z, \quad t\inr
\end{align}
since $h(t) = ct^2$ and $B_2(t) = \N t,\, t\in\R$.
In the following theorem, we consider a general drift function $h$ such that $h(t) +t^2$ is convex on ${E/\R}$, which leads equivalently that $f_\N(t)$ is continuous and convex on $E/\R$. Let thus $s_1 < s_2$ and $t_1 < t_2$  be the two random solutions of $f_\N(s)=0,\, s\in E$ and $f_\N(t)=0, t\inr$, respectively if it holds that
\begin{align*}
f_\N(s^*) = \min_{s\in E}f_\N(s) <0, \quad f_\N(t^*) = \min_{t\inr}f_\N(t) <0.
\end{align*}
\BT\label{T_two}
Let $\MB_\alpha^h(x, E)$ be given by \eqref{MB_E} with $E = [a,b], a<b, a,b\inr$.

(i) If $h(t)+ t^2$ is a continuous, convex function on $E$, then
\begin{align*}
\small \lefteqn{\MB_2^h(x, E) = \int_\R e^z\left[\pk{f_\N(a) < 0, f_\N(b)<0} \td z +\pk{f_\N(a)>0, f_\N(b)>0, s_2 - s_1 >x, f_\N(s^*)<0}\right] \td z} \\
&\quad +\int_\R e^z\left[\pk{f_\N(a)\le0, f_\N(a+x)<0, f_\N(b)>0} + \pk{f_\N(a) > 0, f_\N(b-x)<0, f_\N(b)\le 0}\right] \td z.
\end{align*}
(ii) If $h(t)+t^2, t\inr$ is continuous and convex, and the finite right derivative $h'_+(a)$ and the left derivative $h'_-(b)$ exist with finite values, then
$h'_+(a) \le h'_-(b)$ and
\begin{align*}
\small \lefteqn{\MB_2^h(x, E) =  \int_\R e^z\left[\pk{f_\N(a+x)<0, \sqrt2\N\le h'_+(a) +2a} +\pk{ f_\N(b-x)<0, \sqrt2\N\ge h'_-(b) + 2b}\right] \td z} \\
&\quad + \int_\R e^z\pk{\min(b,t_2) - \max(a,t_1) >x, h'_+(a) +2a <\sqrt2\N < h'_-(b) + 2b, f_\N(t^*)<0} \td z.
\end{align*}
\ET
Typical examples of \netheo{T_two} are discussed in Propositions \ref{T_12} and \ref{T_34} where we take $h(t) = c|t|^\lambda -t^2$ with $\lambda = 1, 2$ and all $\lambda\ge1$. Some \ree{bounds}  are also specified in \neprop{T5}. {All alternative results for $h$ such that  $h(t) + t^2$ is concave are established, see \netheo{Cor1}, Propositions \ref{Cor2} and \ref{Cor3} in Appendix \ref{Appendix}. }

The study on $\MB_\alpha^h$ for fBms $B_\alpha$ with general $\alpha\in(0,1)\cup(1,2)$ is still open since the Slepian's inequality for extremes of processes is not applicable in the sojourn times setting, see e.g., \cite{debicki2015comparison}.

The rest of this paper is organized as follows. Section \ref{Discussion} gives several typical examples   and its bounds as well. \cL{Section \ref{Numerical Study}}  carries out a small scale of numerical studies to illustrate the findings. All proofs are relegated to Section \ref{Proofs}. We present  Appendix \ref{Appendix} for $\MB_2^h$ with concave drift functions.

\section{Further Discussions and Applications}\label{Discussion}
Clearly, a straightforward application of \netheo{T_two} with $h(t) = c|t|^\lambda - t^2, \, t\inr, c>0, \lambda\ge 1$ implies the explicit expressions of $\MB_2^h$, which are given in \neprop{T_12} and \ref{T_34} for  $\lambda =1, 2$ and $\lambda\ge1$. Some bounds are derived in \neprop{T5} for  $\MB_2^h$.
\subsection{Explicit expressions of $\MB_2^h(x,E)$ with $h(t) = c|t|^\lambda -t^2$}
Recalling that $f_\N(t) = h(t)+t^2 -\sqrt2\N t+z = f_{-\N}(-t)$ for symmetry $h$ implies that for the sojourn time $\mathcal L_{\N}(z) = \int_a^{b} \mathbb I(f_\N(t) <0) \td t$
\begin{align*}
\MB_2^h(x, [a,b]) =  \int_{\R} e^{z} \pk{\mathcal L_{\N}(z) > x} \td z = \MB_2^h(x, [-b,-a]),\quad x\in[0, b-a].
\end{align*}
We consider only  $E=[a, b]$ with $b>0$ unless  stated otherwise.
\BP\label{T_12}
Let {$ \MB_2^h(x,E)$} be  the Piterbarg-Berman function  defined in \eqref{MB_E} with $h(t) = c|t|^\lambda - t^2,c>0$ and $E=[a,b]$. Denote by  $f_\N(t) =  c |t|^\lambda - \sqrt 2 \N t + z,\, \N \sim N(0,1)$ as in \eqref{def_f}.

(i) For $\lambda =2$, we have
\begin{align*}
\small \lefteqn{{\MB_2^h(x,E)}= \small \int_\R e^z \left[\pk{f_\N(a+x) <0, \sqrt 2\N/c< {2a+x}} + \pk{f_\N(b-x)<0, \sqrt 2\N/c>{2b-x}}\right]  \td z}\\
& \small + \int_\R e^z \pk{2\N^2/c^2>{ x^2+4 z/c},\,{2a+x}<{\sqrt{2}}\N/c<{2b-x} }  \td z
\qquad\qquad\qquad\qquad
\end{align*}
and for $h\equiv0$,i.e., $c=1$
\begin{align*}
\MB_2(x,[0,T]) = 2\Psi(x/\sqrt2)+ \sqrt2(T-x) \varphi(x/\sqrt2),\quad \MB_2(x) = \lim_{T\to\IF}\frac{\MB_2(x,[0,T])}{T} = \sqrt2\varphi(x/\sqrt2). 
\end{align*}
(ii) For $\lambda=1$, we have with $ \nu(m,c) =  e^{m^2 -c|m|} \Psi\left([c-2m]/{\sqrt2}\right)$\begin{align*}
\MB_2^h(x,{E}) = \small \left\{
\begin{array}{ll}
 e^{(a+x)(a+x-c)}  - \nu(a+x,c) + \nu(b-x,c), & a\ge0,\\
 \nu(b-x, c) + \nu(-(a+x), c) +  \int_{-\IF}^0 e^z\pk{\sqrt2|\N|<c,  \min(b, t_2) - \max(a, t_1) >x} \td z,
 & a<0,
  \end{array}
  \right.
\end{align*}
{where $t_1<t_2$ are as in \netheo{T_two} (ii), i.e., the random solutions of $f_\N(t)=0,\, t\inr$ equal}
\begin{align*}
t_1 = \frac z{c+\sqrt2\N}<0< \frac{-z}{c-\sqrt2\N} = t_2,\quad \cL{z<0}.
\end{align*}
\EP
{Intuitively}, the three parts involved for $\lambda =2$ above are obtained via the comparison among the symmetric axis $t=\sqrt2\N/(2c)$ of the quadratic symmetric curve $f_\N(\cdot)$,  $a + x/2$ and $b-x/2$. Meanwhile,  the well-known Pickands' constant $\MB_2(0)= \sqrt\pi$ can be implied by  $\MB_2(x)= \sqrt2\varphi(x), x\ge0$. The specification for $\lambda=1$ follows from \netheo{T_two} {(ii)}, and the expression of $\MB_2^h(x, E)$ with $E=[a,b]$ including the original point is more involved due to the piece-wise property of $f_\N(\cdot)$. Its detailed expansions are given  in Appendix \ref{append_2}.
Below, we consider the general drift function $h(t) = c|t|^\lambda -t^2$ with $\lambda \ge1$. The case with $0<\lambda<1$ is given in \neprop{Cor2} in Appendix \ref{append_1}.
\BP\label{T_34}
Let $\MB_2^h(x,E)$ be  the Piterbarg-Berman function  defined in \eqref{MB_E} with drift function $h(t)= c|t|^\lambda - t^2, c>0, \lambda \ge 1$ \ree{and $E=[a,b]$}.

(i) For $a\ge0$, we have
\begin{align*}
\lefteqn{ \MB_{2}^h(x,E)
= \E{\mathbb I(\sqrt 2\N (a+x)  + \W> c(a+x)^{\lambda})}}  \\
&\quad +\int_0^{\IF} e^z\pk{\min(b, t_2) -  \max(a, t_1) >x, f_\N(t^*) <0, \N>0}   \td z.
\end{align*}
(ii) For $a<0$, we have
\begin{align*}
\small \lefteqn{\MB_{2}^h(x,E)
= \int_\R e^z\left[\pk{f_\N(a+x)<0, f_\N'(a) {\ge}0}  + \pk{ f_\N(b-x)<0, f_\N'(b) {\le }0} \right] \td z} \\
&\quad +\int_\R e^z\pk{\min(b,t_2) - \max(a,t_1) >x,  f_\N'(a) <0 < f_\N'(b), f_\N(t^*)<0} \td z.
\end{align*}
\cL{Here $t_1< t_2$  are defined as in \netheo{T_two} (ii), i.e., the two random solutions of $f_\N(t)=0, t\inr$ when its minimum $f_\N(t^*)$ is less than zero.}
\EP
We see that the Piterbarg-Berman functions  require  cumbersome calculations even for typical drift functions, see Theorem \ref{T_two}, Propositions \ref{T_12} and \ref{T_34}. Below, we consider alternatively its bounds.
\subsection{Bounds of the Piterbarg-Berman functions} \label{Bounds}
Recalling the geometry property of the convex curve $f_\N$, we develop below an upper bound of $\MB_2^h(x, E)$. To simplify the notation, we consider the setting of \neprop{T_34}, i.e., $h(t) = c|t|^\lambda - t^2$ with $c>0, \lambda\ge1$ {(see {\neprop{Cor3} for the lower bound of $\MB_2^h(x, E)$ with $0< \lambda <1$} in Appendix \ref{append_1}).} The general bounds $\MB_\alpha^h(x, E)$ for $\alpha\in(0,2]$ are also established by linking the degenerate cases with $h, x, E$ taken into consideration.
Set with $c_0 = c_0(y)  = \lambda c \abs{y}^{\lambda-1}\sign(y)$ and $\nu'(m, c) = e^{m^2- c m }\Psi([c - 2m]/\sqrt2)$
\begin{align}\label{Eq_C0}
\qquad C_0(y) = \expon{(\lambda -1) c\abs{y}^{\lambda}}, \quad D_0(y, E)  =  \nu'(b-x, c_0) + \nu'(-(a+x), -c_0), \quad {E=[a,b]}.
\end{align}
\BP
\label{T5}
Let $\MB_\alpha^h(x, E)$
be the Piterbarg-Berman function given by \eqref{MB_E}.

(i) For $\alpha=2$ and $h(t) = c|t|^\lambda - t^2, c>0, \lambda \ge1$, we have for $E=[a,b]$
\begin{align*}
\MB_2^h(x,{E})  \le   \min_{y\in{E}\setminus\{0\}} C_0(y) D_0(y, {E}).
\end{align*}
(ii) For $\alpha \in (0,2]$ and drift function $h$ satisfying $M= \max_{s\in E}h(t) <\IF, E \subset \R_+\cup\{0\}$. We have
\begin{align*}
  e^{-M} \MB_\alpha(x, E)  \le \MB_\alpha^h(x,E)
  \le \min(\MB_\alpha^h(x,[0,\IF]), \MB_\alpha^h(0,E)), \quad x\ge0.
\end{align*}
\EP
\vspace{-0.82cm}
\section{Numerical study}\label{Numerical Study}
In order to illustrate the theoretical findings in Theorem \ref{T_one_two}, Propositions \ref{T_12} and \ref{T5}, we carry out a small scale of numerical studies for the Piterbarg-Berman function $\MB_\alpha^h(x, E)$.
\COM{We summarize it as follows.
\begin{enumerate}
\item[(i)] In  Figure \ref{fig_Berman_alpha_1}, we consider $E=[0,\IF]$ for $\MB_1^h(x, E)$. Applying  \netheo{T_one_two} for $h(t) = c|t|$ with different $c$'s, we see that the larger the $c$ is, the more quickly the curve decreases with respect to $x$ (the same as below), and these curves with different $c$'s become closer and closer as $x\to0$ and $x\to\IF$.
\item[(ii)]  In Figure \ref{fig_Berman_two} (and thereafter), we consider finite time interval $E=[a,b]$ and $h(t) = c|t|^\lambda - t^2,\, c>0$ as in \neprop{T_12}.  Here we take $\lambda =2, E = [0,3], [-2,1], [1,3], [-1,1],[a,b]= [0,4]$ and  $c=1, 2, 4$.  We see that $\MB_2^h(x, E), x\in[0,b-a]$ is decreasing slowly.
\item[(iii)]  In Figure \ref{fig_Berman_one},  we draw $\MB_2^h(x,E)$ with $h(t) = c|t| -t^2$. It seems more sensitive to the time interval $E$ of the same length and decreases strongly as $x\to0$ but indifferent for larger $x$.
\item[(iv)] In Figure \ref{fig_Upper_Bound}, we apply \neprop{T5} and show the efficiency of the upper bounds of $\MB_2^h(x, E)$ with $h(t) = (c-1)t^2$ for $c=1, 2$. Clearly, the upper bound of the Piterbarg-Berman {function} become closer to the true values  for larger $x$.  The error is also determined by the time intervals and the coefficient $c$ of the drift function as well.
\end{enumerate}}



\begin{figure}[H]
\vspace{-1cm}
	\centering
	\subfigure{
		\begin{minipage}[t]{0.5\linewidth}
			\centering
			\includegraphics[width=8cm,height=6cm]
			{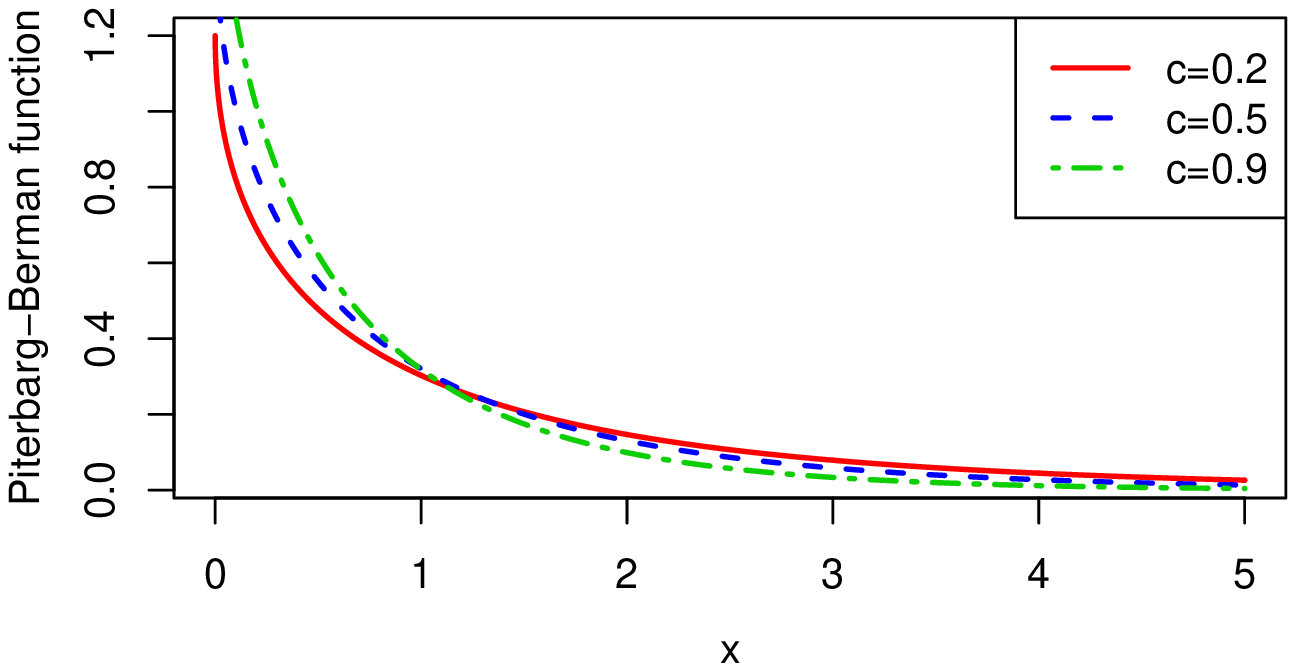}
			\caption*{}
		\end{minipage}%
	}%
			\subfigure{
		\begin{minipage}[t]{0.5\linewidth}
			\centering
			\includegraphics[width=8cm,height=6cm]{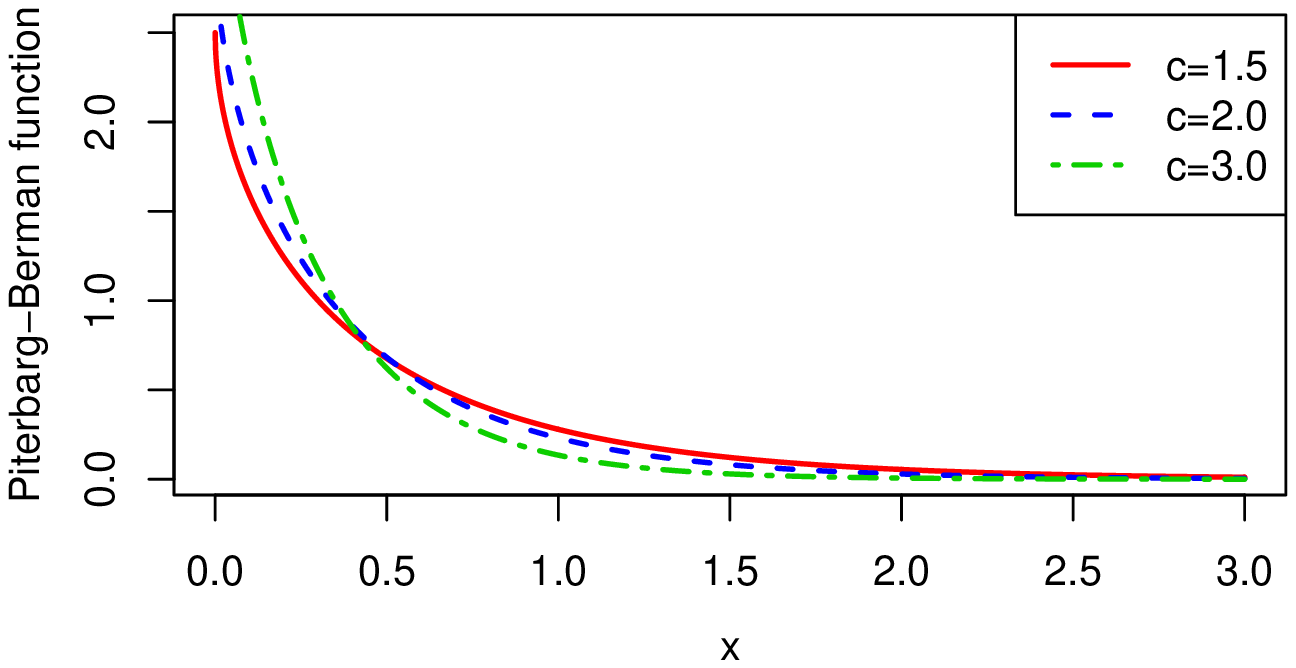}
			\caption*{}
		\end{minipage}%
	}%
\vspace{-2cm}
	\centering
	\caption{\cL{The Piterbarg-Berman function {$\MB_{1}^h(x)$}
			with $h(t) = ct$ for $c\in(0,1)$ (left) and $c\in (1,\IF)$ (right).} }
	\label{fig_Berman_alpha_1}
\vspace{-1cm}
\end{figure}
In  Figure \ref{fig_Berman_alpha_1}, we consider $E=[0,\IF]$ for $\MB_1^h(x, E)$. Applying  \netheo{T_one_two} for $h(t) = c|t|$ with different $c$'s, we see that the larger the $c$ is, the more quickly the curve decreases with respect to $x$ (the same as below), and these curves with different $c$'s become closer and closer as $x\to0$ and $x\to\IF$.
\begin{figure}[H]
\vspace{-1cm}
	\centering
	\subfigure{
		\begin{minipage}[t]{0.5\linewidth}
			\centering
			\includegraphics[width=8cm,height=6cm]{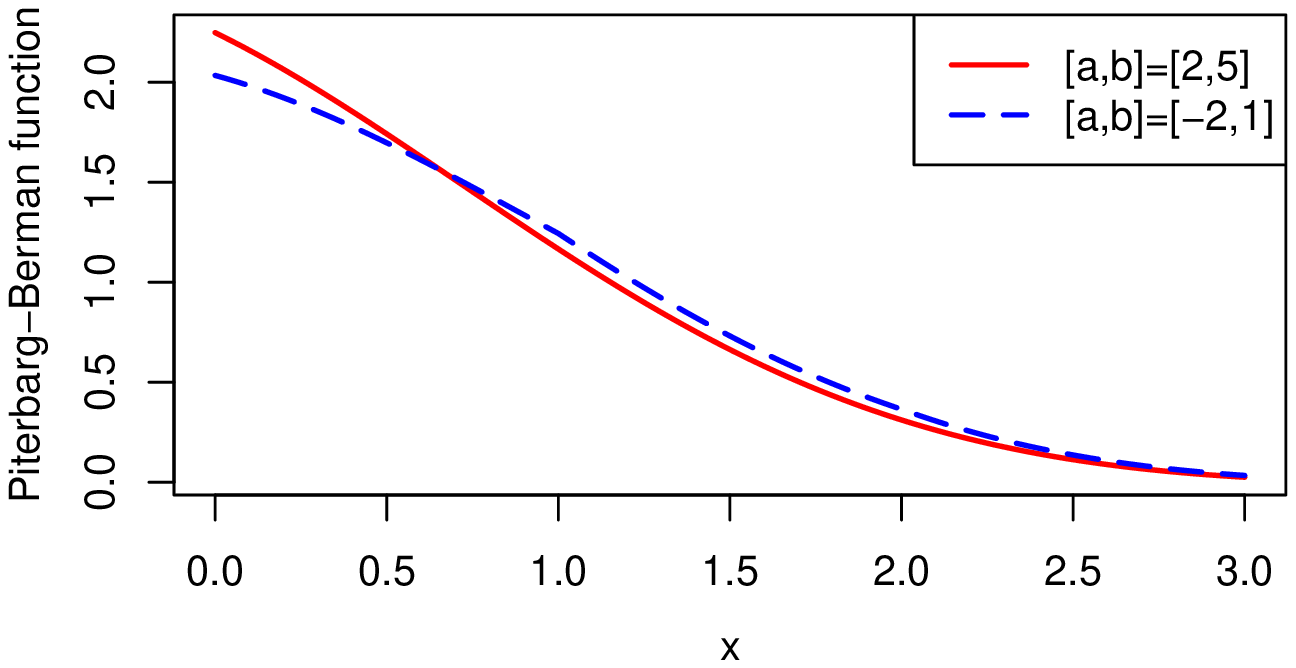}
			\caption*{}
		\end{minipage}%
	}%
	\subfigure{
		\begin{minipage}[t]{0.5\linewidth}
			\centering
			\includegraphics[width=8cm,height=6cm]{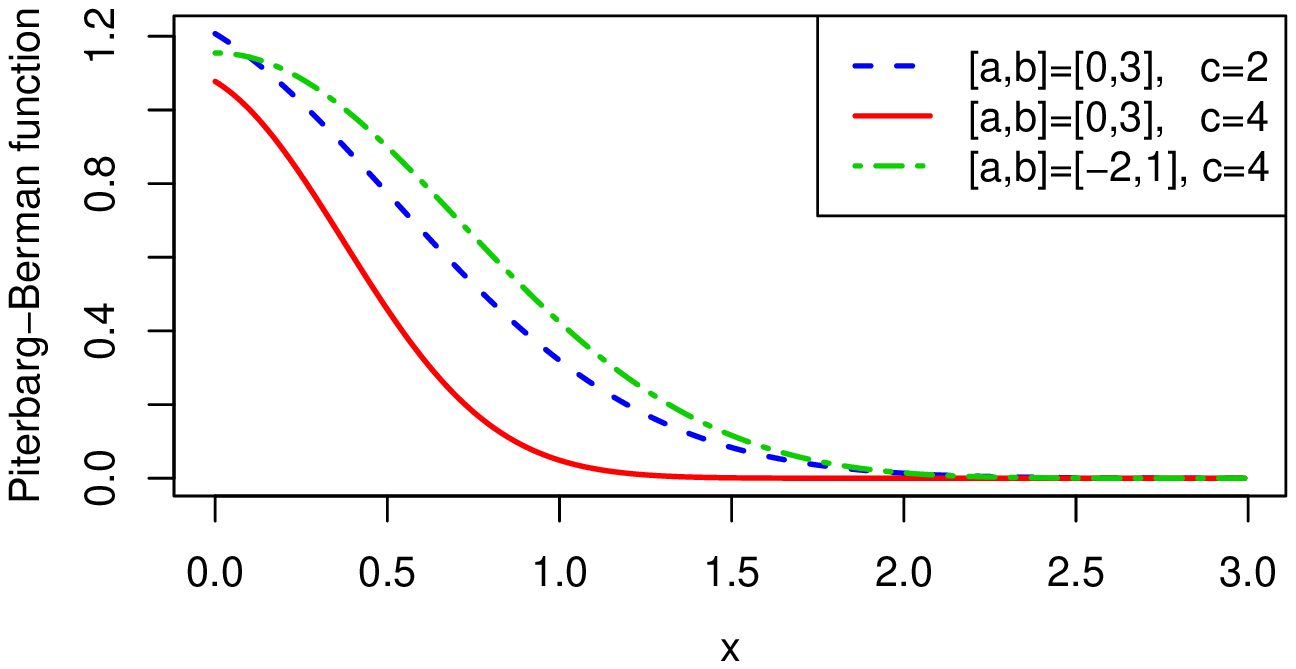}
			\caption*{}
		\end{minipage}%
	}%
\vspace{-2cm}
	\centering
	\caption{\cL{The Piterbarg-Berman function {$\MB_{2}^h(x,{E})$} 
			with $h(t) = c \abs t^\lambda -t^2$ for $\lambda =2, c=1$(left) and $\lambda =2, c=2,4$ (right).} }  
	\label{fig_Berman_two}   
\vspace{-1cm}
\end{figure}
In Figure \ref{fig_Berman_two} (and thereafter), we consider finite time interval $E=[a,b]$ and $h(t) = c|t|^\lambda - t^2,\, c>0$ as in \neprop{T_12}.  Here we take $\lambda =2, E = [0,3], [-2,1], [1,3], [-1,1], [0,4]$ and  $c=1, 2, 4$.  We see that $\MB_2^h(x, E), x\in[0,b-a]$ is decreasing slowly.
\begin{figure}[H]
\vspace{-0.5cm}
	\centering
	\subfigure{
		\begin{minipage}[t]{0.5\linewidth}
			\centering
			\includegraphics[width=8cm,height=6cm]{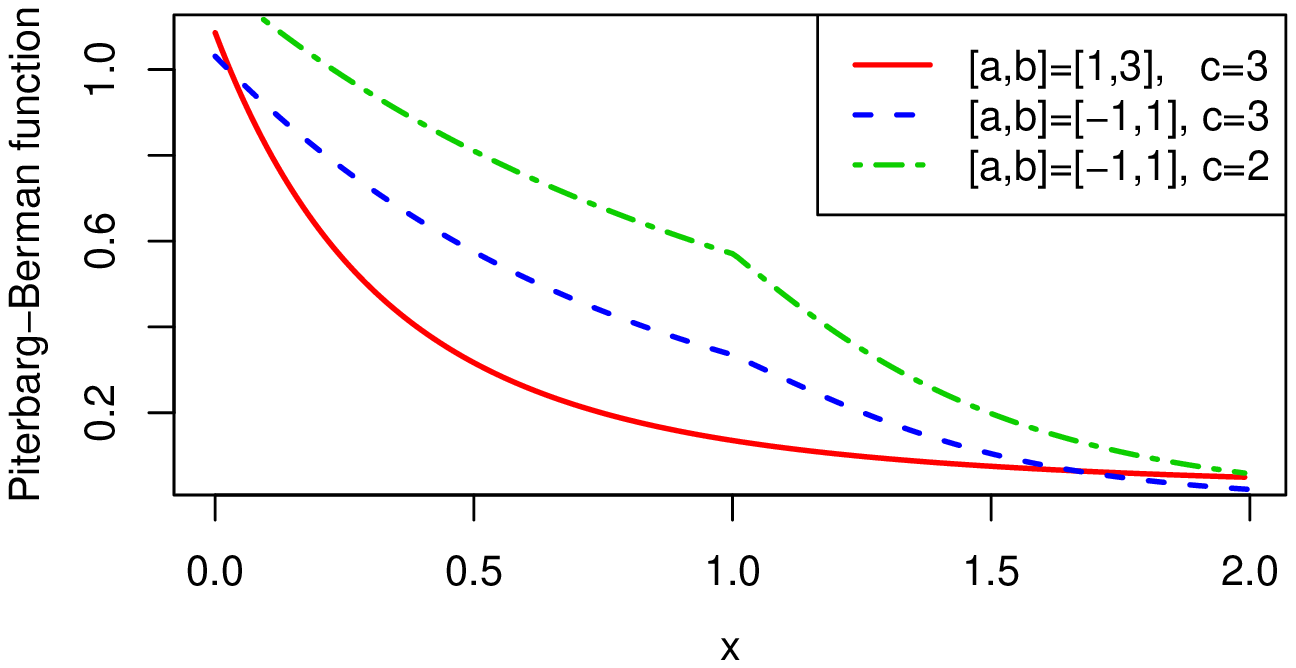}
			\caption*{}
		\end{minipage}%
	}%
	\subfigure{
		\begin{minipage}[t]{0.5\linewidth}
			\centering
			\includegraphics[width=8cm,height=6cm]{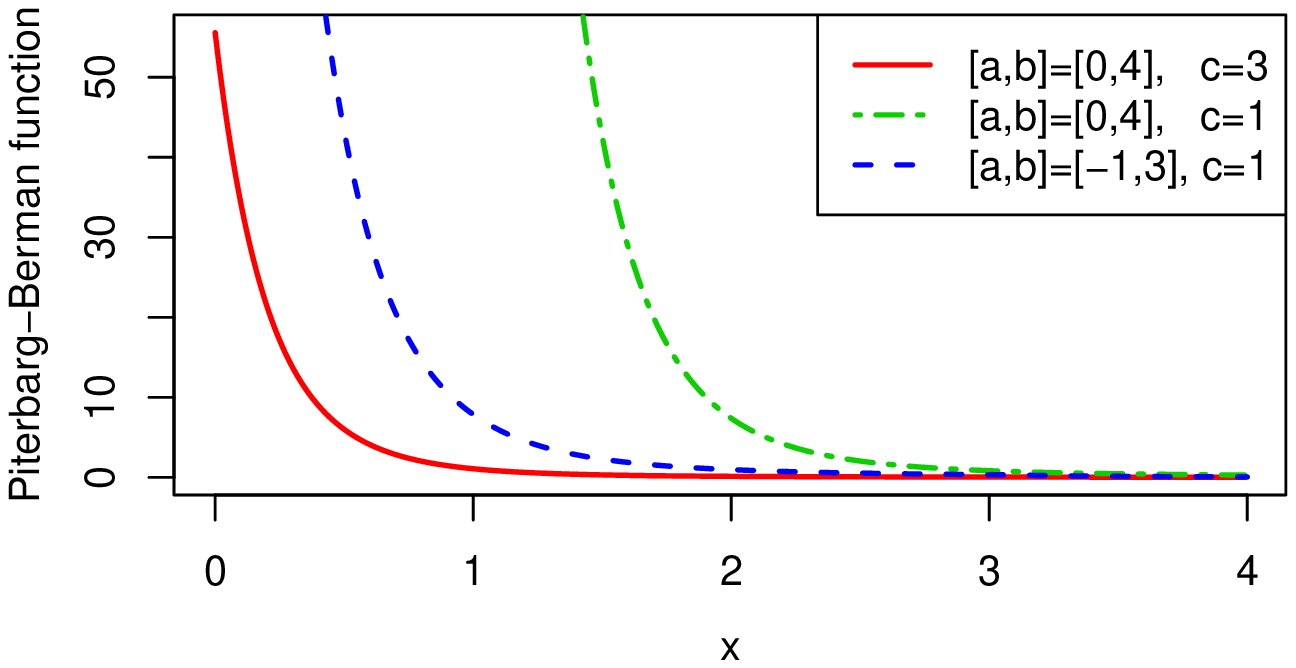}
	\end{minipage}%
	}%
\vspace{-2cm}
	\centering
	\caption{\cL{The Piterbarg-Berman function  $\MB_{2}^h(x,{E})$ with drift function $h(t) = c \abs t-t^2$.}}
	\label{fig_Berman_one}
\vspace{-0.8cm}
\end{figure}
In Figure \ref{fig_Berman_one}, we draw $\MB_2^h(x,E), \ree{E=[a,b]}$ with $h(t) = c|t| -t^2$. It seems more sensitive to the time interval $E$ of the same length and decreases strongly as $x\to0$ but indifferent for larger $x$.
\begin{figure}[H]
\vspace{-0.6cm}
	\centering
	\subfigure{
		\begin{minipage}[t]{0.5\linewidth}
			\centering
			\includegraphics[width=8cm,height=5cm]{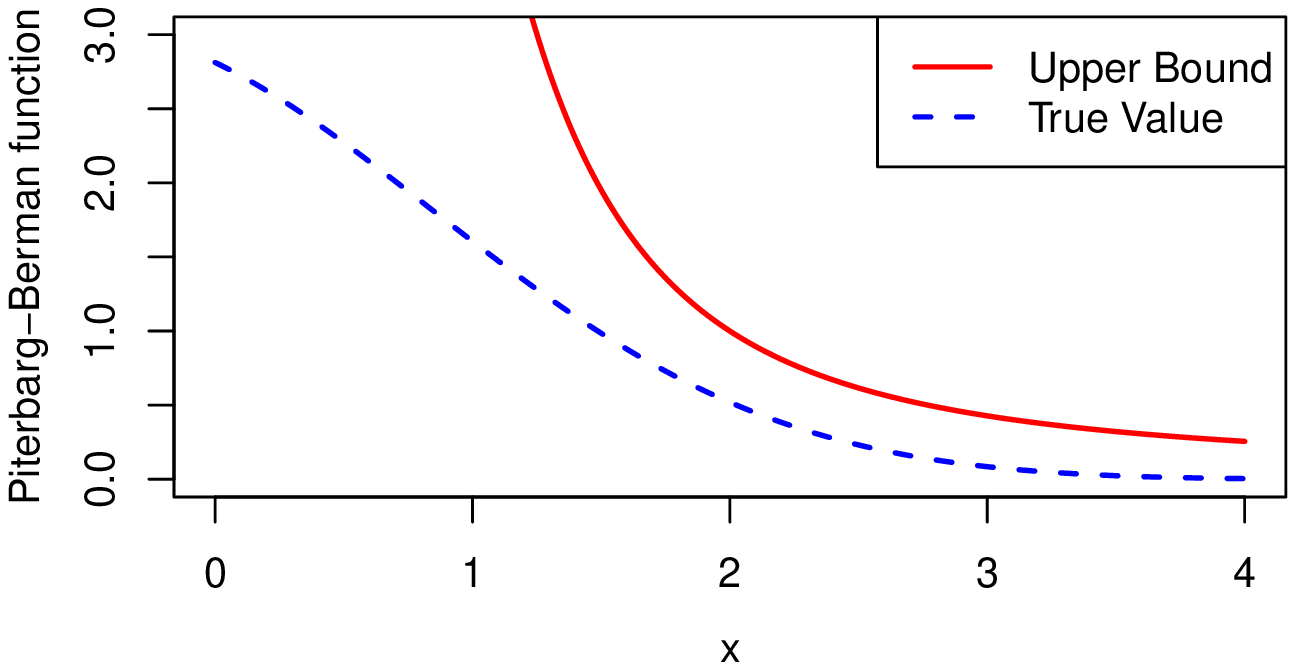}
			\vspace{-0.7cm}
			\caption*{(i) $[a, b]=[1,5]$}
		\end{minipage}%
	}%
	\subfigure{
		\begin{minipage}[t]{0.5\linewidth}
			\centering
			\includegraphics[width=8cm,height=5cm]{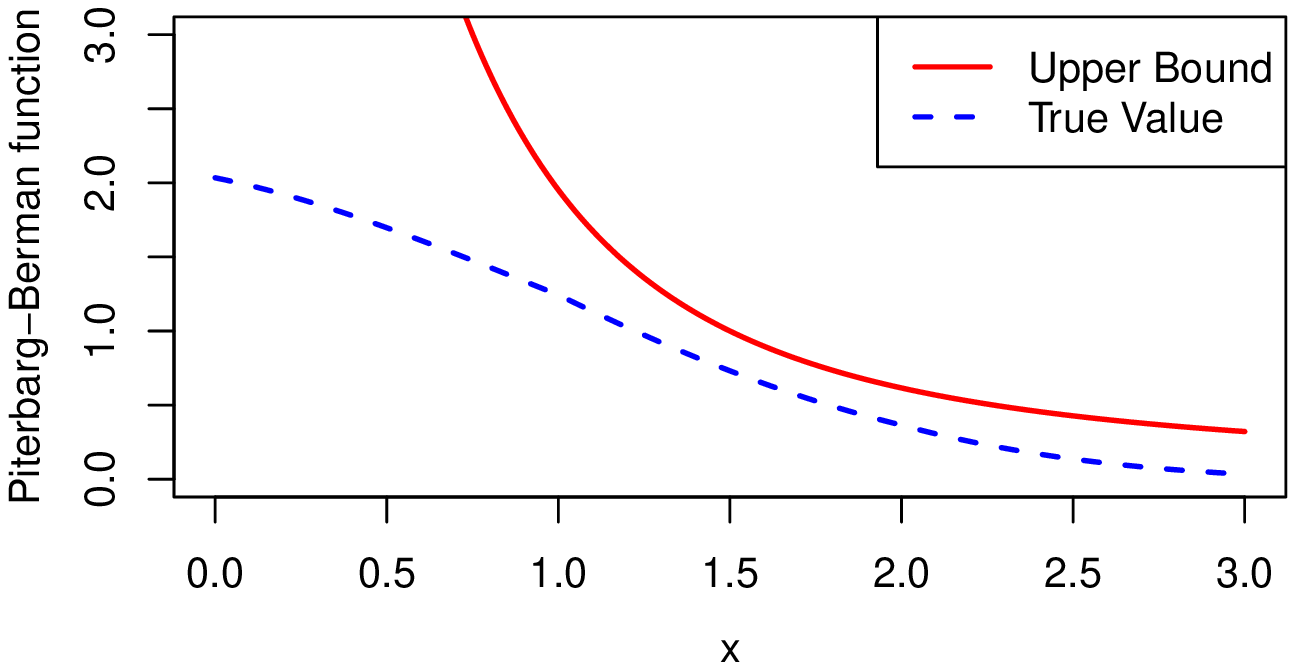}
			\vspace{-0.7cm}
			\caption*{(ii) $[a,b] =[-1,2]$}
		\end{minipage}%
	}\\
\vspace{-0.5cm}
\subfigure{
	\begin{minipage}[t]{0.5\linewidth}
		\centering		\includegraphics[width=8cm,height=5cm]{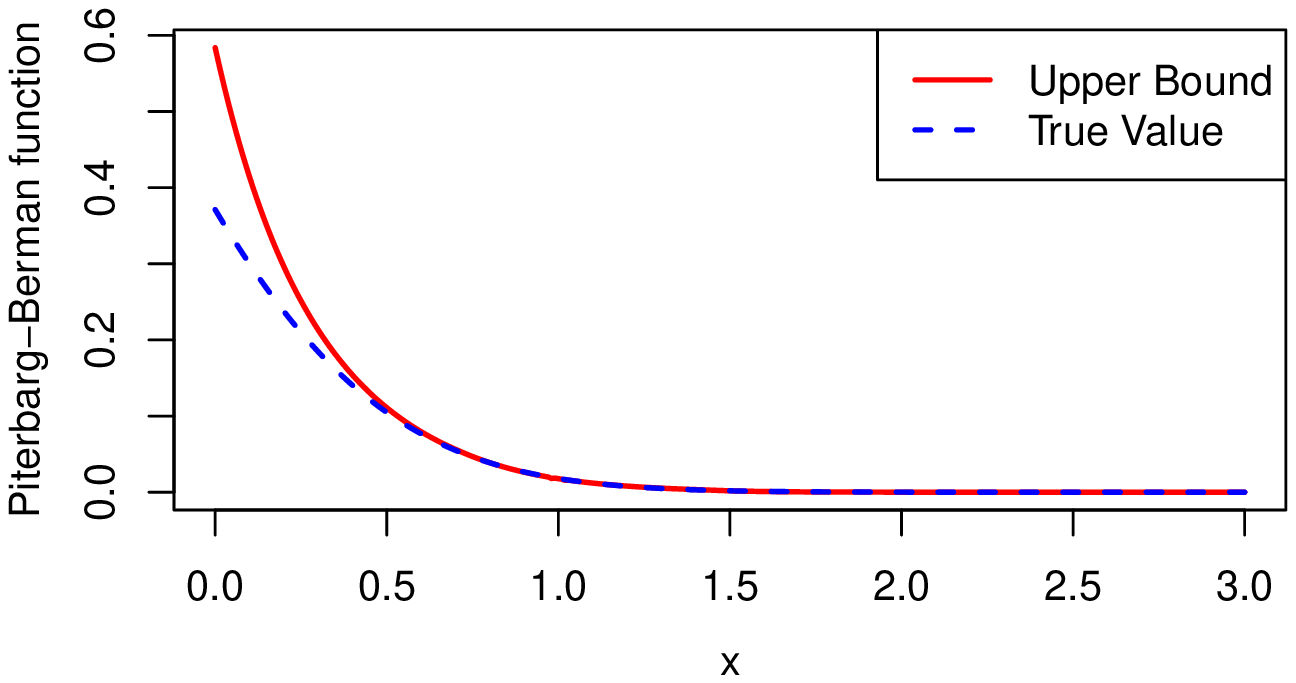}
		\vspace{-0.7cm}
		\caption*{(iii) $ [a, b]=[1, 4]$}
	\end{minipage}%
}%
\subfigure {
	\begin{minipage}[t]{0.5\linewidth}
		\centering		\includegraphics[width=8cm,height=5cm]{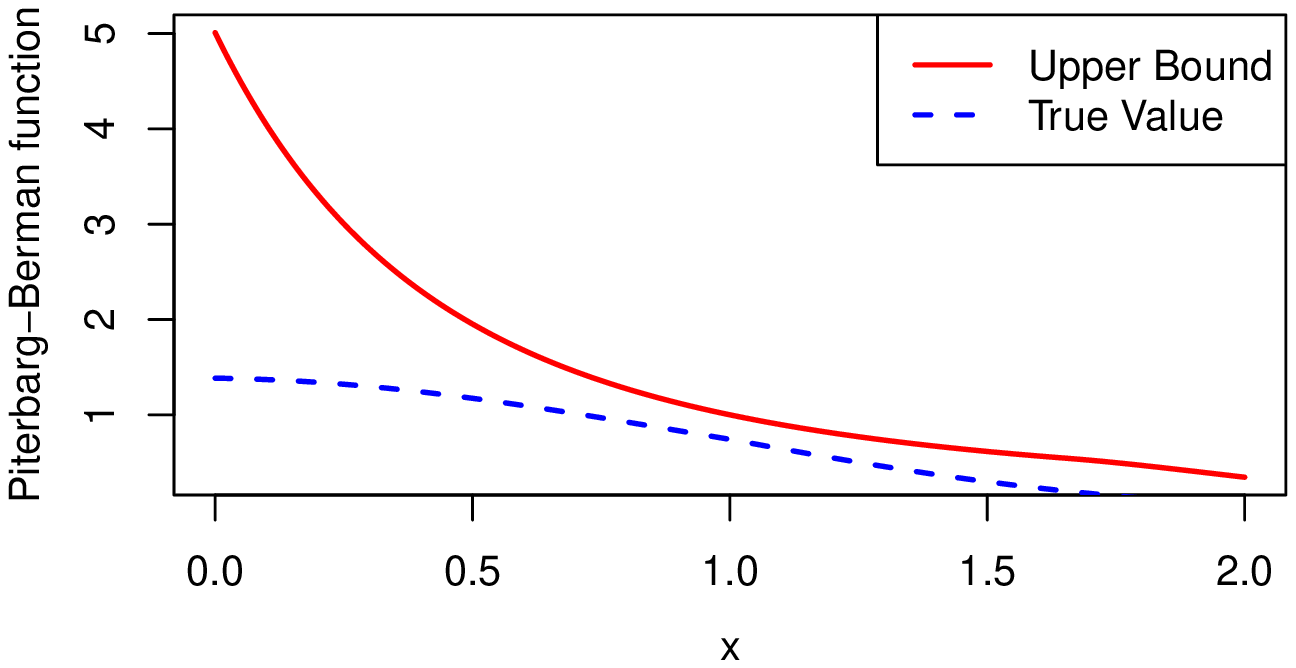}
		\vspace{-0.7cm}
		\caption*{(iv) $[a,b]=[-1, 1]$}
	\end{minipage}%
}%
\vspace{-0.5cm}
\centering
	\caption{Upper bounds and true values of $\MB_2^h(x, {E})$  for drift function $h(t) = c|t|^\lambda - t^2$ with $\lambda=2, c=1$ (top) and $\lambda=2, c=2$ (bottom).}\label{fig_Upper_Bound}
\vspace{-0.7cm}
\end{figure}

In Figure \ref{fig_Upper_Bound}, we apply \neprop{T5} and show the efficiency of the upper bounds of $\MB_2^h(x, E), \ree{E=[a,b]}$ with $h(t) = (c-1)t^2$ for $c=1, 2$. Clearly, the upper bound of the Piterbarg-Berman {function} become closer to the true values  for larger $x$.  The error is also determined by the time intervals and the coefficient $c$ of the drift function as well.

Overall, the Piterbarg-Berman {function} has nice properties with respect to the argument $x$ and the coefficient $c$ arising in the simple drift {function} $h(t)= c|t|^\lambda - t^2$. Meanwhile, its essential complexity is closely related to the observed time interval ${E=}[a,b]$ as well as the power index $\lambda$ involved.

\section{Proofs}\label{Proofs}
Throughout the proofs, we keep the same notation aforementioned unless stated otherwise. We shall present the proofs of Theorem \ref{T_one_two} and Propositions \ref{T_12}-\ref{T5} one by one.
\subsection{Proofs of \netheo{T_one_two} for $\MB_\alpha^h(x)$ with $h(t) = ct^\alpha, c>0$ and $\alpha =1, 2$}
\ \newline
\underline{\textbf{Proof of (i) $\alpha=1$}}. Recalling the sojourn time $Y_z$ given by \eqref{def_sojourn_BM}, we have
\begin{align}\label{decomp}
\MB_1^h(x) = \left(\int_0^\IF + \int_{-\IF}^0\right) e^z\pk{Y_z > x} \td z =: I_1(x) + I_2(x).
\end{align}
We deal with the integrals $I_1(x)$ and $I_{2}(x)$ according to the upward and downward crossings.

{For $I_1(x)$}. By the lack of upward jumps and the strong Markov property, we have for $z\ge0$
\begin{align*}
\pk{Y_z > x} = \pk{Y_0 >x} \pk{\sup_{t\ge0} (\sqrt2 B_1(t)-\cL{(1+c)t}) \ge z},
\end{align*}
where, it follows from \cite{Borodin2012} that (see e.g., Eq. (3) in p. 255 therein)
\begin{align}
\label{def_survival_Y0}
&\pk{\sup_{t\ge0} (\sqrt2B_1(t)-\cL{(1+c)}t) \ge z} = e^{-\cL{(1+c)}z},\quad z\ge0, \notag\\
&\pk{Y_0 >x} = 2\left(1+\frac{(1+c)^2x}2\right) \Psi\left((1+c)\sqrt{\frac x2}\right) - (1+c)\sqrt{2x} \varphi\left((1+c)\sqrt{x/2}\right)\notag
\\
&\qquad\qquad\quad = 2(1+x_{c}^2)\Psi(x_{\cL{c}})-2x_{\cL{c}}\varphi(x_{\cL{c}}),\quad x_c = (1+c)\sqrt{x/2}.
\end{align}
Therefore,
\begin{align}\label{Ix}
I_1(x) = \int_{0}^\IF e^{-\cL{c}z}\td z \pk{Y_0>x} =  \frac2{\cL{c}}\left[ (1+x_{c}^2)\Psi(x_{\cL{c}})-x_{\cL{c}}\varphi(x_{\cL{c}})\right].
\end{align}
{For $I_{2}(x)$}. Note that  $\tau_z = \inf\{t\ge0: \sqrt2 B_1(t)-\cL{(1+c)}t \ge z\}$ for given $z<0$  is a stopping time with cumulative distribution function (cdf) and probability  density function (pdf) given by
\begin{align}
\label{def_survival_tau}
\begin{split}
\pk{\tau_z \le t} &= \Psi\left(\frac{-z}{\sqrt{2t}} - t_{\cL{c}}\right)  + e^{-\cL{(1+c)}z} \Psi\left(\frac{-z}{\sqrt{2t}} + t_{\cL{c}}\right), \quad z<0,\\
\cL{f_z(t)} &= \frac{-z}{\sqrt{2 t^3}} \varphi\left(\frac{-z}{\sqrt{2t}} - t_{\cL{c}}\right), \quad t_c = (1+c)\sqrt{x/2}, \, \, t >0
\end{split}
\end{align}
by Eq. (3) in page 261 of \cite{Borodin2012} and the fact that $ \varphi(z/\sqrt{2t} + t_{\cL{c}}) = e^{\cL{(1+c)}z} \varphi(z/\sqrt{2t} - t_{\cL{c}})$. Therefore, again by the strong Markov property, we have
\begin{align*}
\pk{Y_z > x} &= \pk{ \tau_z + \int_{\cL{\tau_z}}^\IF\mathbb I(\sqrt2 B_1(t)\cL{-(1+c)t} > z) \td t > x} \\
&= \pk{ \tau_z + \int_0^\IF \mathbb I(\sqrt2 [B_1(t+\tau_z) - B_1(\tau_z)]\cL{-(1+c)t} > 0) \td t > x} \\
&= \pk{\tau_z > x} + \int_0^x \pk{Y_0 > x-t} \cL{f_z(t)}\td t,
\end{align*}
where the last step follows since $\tau_z$ and $Y_0$ are independent.
Consequently,
\begin{align}
\label{def_II2}
I_{2}(x) &= \int_{0}^\IF e^{-z} \left[\pk{\tau_{-z} > x} + \int_0^x \pk{Y_0 > x-t} \cL{f_{-z}(t)}\td t\right] \td z
\notag\\
&=: I_{21}(x) + I_{22}(x).
\end{align}

In the following, we deal with $I_{21}(x)$ and $I_{22}(x)$ subsequently. First, we have by \eqref{def_survival_tau}
\begin{align}\label{II1}
I_{21}(x) &= 1- \int_0^\IF e^{-z}\left[\Psi\left(\frac{z}{\sqrt{2x}} - x_{\cL{c}}\right)  + \cL{e^{(1+c)z}} \Psi\left(\frac{z}{\sqrt{2x}} + x_{\cL{c}}\right)\right] \td z \notag\\
\COM{&= 1- \int_0^\IF \varphi(u-x_{\cL{c}})\int_0^{\sqrt{2x}u} e^{-z}\td z - \int_0^\IF \varphi(u+x_{\cL{c}})\int_0^{\sqrt{2x}u} e^{\cL{c}z}\td z \notag\\
&= 1-\int_0^\IF[1-e^{-\sqrt{2x}u}]\varphi(u-x_{\cL{c}}) \td u + \frac1{\cL{c}}\int_0^\IF [1-e^{\cL{c}\sqrt{2x}u}]\varphi(u+x_{\cL{c}}) \td u \notag\\
&= 1-\int_0^\IF[\varphi(u-x_{\cL{c}})-\varphi(u+x_{\cL{c}}')e^{-\cL{c}x}] \td u + \frac1{\cL{c}}\int_0^\IF [\varphi(u+x_{\cL{c}})-\varphi(u+x_{\cL{c}}')e^{-\cL{c}x}]\td u \notag\\
&= 1- \Psi(-x_{\cL{c}})+\frac1{\cL{c}}\Psi(x_{\cL{c}}) + \left[1-\frac1{\cL{c}}\right]\Psi(x_{c}')e^{-\cL{c}x} \notag\\}
&= \frac{\cL{1+c}}{\cL{c}}\Psi(x_{\cL{c}}) - \frac {\cL{1-c}}{\cL{c}}\Psi(x_{c}')e^{-\cL{c}x}, \quad x_c' = (1-c)\sqrt{x/2}.
 \end{align}
 Similarly, we have by \eqref{def_survival_Y0}, \eqref{def_survival_tau} and \eqref{def_II2}
\begin{align*}
I_{22}(x) &=  \int_0^x \pk{Y_0 > x-t} \Big[\int_0^\IF e^{-z}\cL{f_{-z}(t)} \td z\Big] \td t \\
&=  2\int_0^x \left\{[(1+x_{c}^2)-t_{\cL{c}}^2]\Psi(\sqrt{x_{c}^2-t_{\cL{c}}^2})-\sqrt{x_{c}^2-t_{\cL{c}}^2}\varphi(\sqrt{x_{c}^2-t_{\cL{c}}^2})\right\}\Big[\cL{(1+c)}\frac{\varphi(t_{\cL{c}})}{t_{\cL{c}}} - \cL{(1-c)} \Psi(t_{c}')e^{-\cL{c}t}\Big]\td t\\
\COM{&= 2\int_0^x\left\{\cL{(1+c)}\left[ (1+x_{c}^2)\frac{\varphi(t_{\cL{c}})}{t_{\cL{c}}}\Psi(\sqrt{x_{c}^2-t_{\cL{c}}^2}) - t_{\cL{c}} \varphi(t_{\cL{c}})\Psi(\sqrt{x_{c}^2-t_{\cL{c}}^2}) - \sqrt{\frac{x}{t}-1} \varphi(t_{\cL{c}})\varphi(\sqrt{x_{c}^2-t_{\cL{c}}^2})\right]\right.\\
&\quad -\cL{(1-c)}\left[(1+x_{c}^2)\Psi(\sqrt{x_{c}^2-t_{\cL{c}}^2})\Psi(t_{c}')e^{-\cL{c}t} -t_{\cL{c}}^2 \Psi(\sqrt{x_{c}^2-t_{\cL{c}}^2})\Psi(t_{c}')e^{-\cL{c}t}\right.\\
&\quad -\left.\left. \sqrt{x_{c}^2-t_{\cL{c}}^2}\varphi(\sqrt{x_{c}^2-t_{\cL{c}}^2})\Psi(t_{c}')e^{-\cL{c}t}\right]\right\}\td t \\}
&=: 2\Big\{\cL{(1+c)}\big[(1+x_{c}^2)A_1 - A_2 - A_3\big] - \cL{|1-c|}\I{c\neq1}\big[(1+x_{c}^2)A_4 - A_5 - A_6\big]\\
 &\quad + \cL{|1-c|}\I{\cL{c>1}} \big[(1+x_{c}^2)A_4^+ - A_5^+ - A_6^+\big]\Big\},
 \end{align*}
where, with $t_c = (1+c)\sqrt{t/2}, t_c'=(1-c)\sqrt{t/2}$,
\begin{align*}
A_1&=\int_0^x\frac{\varphi(t_{\cL{c}})}{t_{\cL{c}}}\Psi(\sqrt{x_{c}^2-t_{\cL{c}}^2})\td t,\quad A_2 = \int_0^x t_{\cL{c}} {\varphi(t_{\cL{c}})}\Psi(\sqrt{x_{c}^2-t_{\cL{c}}^2})\td t,\quad A_3 = \frac{e^{-x_{c}^2/2}}{2\pi}\int_0^x \sqrt{\frac{x}{t}-1} \td t,\\
A_4&=\int_0^x \Psi(\sqrt{x_{c}^2-t_{\cL{c}}^2})\Psi(|t_{c}'|)e^{-\cL{c}t} \td t,\qquad\qquad\qquad A_4^+=\int_0^x \Psi(\sqrt{x_{c}^2-t_{\cL{c}}^2})e^{-\cL{c}t} \td t, \\
A_5&=\int_0^x t_{\cL{c}}^2\Psi(\sqrt{x_{c}^2-t_{\cL{c}}^2})\Psi(|t_{c}'|)e^{-\cL{c}t} \td t,\quad\qquad\qquad A_5^+=\int_0^x t_{\cL{c}}^2\Psi(\sqrt{x_{c}^2-t_{\cL{c}}^2})e^{-\cL{c}t} \td t, \\
A_6&=\int_0^x \sqrt{x_{c}^2-t_{\cL{c}}^2}\varphi(\sqrt{x_{c}^2-t_{\cL{c}}^2})\Psi(|t_{c}'|)e^{-\cL{c}t} \td t,\qquad A_6^+=\int_0^x \sqrt{x_{c}^2-t_{\cL{c}}^2}\varphi(\sqrt{x_{c}^2-t_{\cL{c}}^2})e^{-\cL{c}t} \td t.
\end{align*}
Here the three terms $A_4^+ \sim A_6^+$ for $c>1$ arises since $\Psi(t_c') = 1- \Psi(|t_c'|)$ as $t_c'<0$, i.e., $c>1$.

Next,  we deal with $I_{22}(x)$ by specifying the three sets of integrals: (i) $A_1\sim A_3$ with $c>0$; (ii) $A_4\sim A_6$ with $c\neq1$ and (iii) $A_4^+\sim A_6^+$ with  $c>1$.

\underline{(ii) For $A_1\sim A_3$ with $c>0$}. By the symmetry of standard normal distributions, we have
\begin{align}
\cL{(1+c)^2}A_1 &=  4\int_0^{x_{\cL{c}}}\varphi(v)\Psi(\sqrt{x_{c}^2-v^2})\td v \notag\\
&= \left(\int_{u^2+v^2 > x_{c}^2} - \int_{v^2 > x_{c}^2}\right)\varphi(u)\varphi(v) \td u\td v = e^{-x_{c}^2/2}-2\Psi(x_{\cL{c}}),\\
\label{A2}
\cL{(1+c)^2}A_2 &= 4\int_{u^2+v^2 > x_{c}^2, |v|<x_{\cL{c}}} v^2\varphi(v)\varphi(u) \td u \td v  = \left[\int_{u^2+v^2>x_{c}^2} - \int_{|v|>x_{\cL{c}}}\right]v^2\varphi(v)\varphi(u) \td u \td v \notag\\
&= \frac12 \int_{u^2+v^2>x_{c}^2} (u^2+v^2)\varphi(u)\varphi(v) \td u \td v - 2[\Phi(v)-v\varphi(v)]|_{x_{\cL{c}}}^\IF\notag \\
&= -(1+t)e^{-t}|_{x_{c}^2/2}^\IF - 2[\Phi(v)-v\varphi(v)]|_{x_{\cL{c}}}^\IF \notag\\
&= [1+x_{c}^2/2]e^{-x_{c}^2/2} - 2[\Psi(x_{\cL{c}})+ x_{\cL{c}}\varphi(x_{\cL{c}})]=: (-2)\mathrm m(x_c).
\end{align}
For $A_3$, we take $s=\sqrt{t/x}, 0 < t <x$ and thus
\begin{align*}
A_3 = \frac{x}{\pi}e^{-x_{c}^2/2}\int_0^1 \sqrt{1-s^2}\td t = \frac{xe^{-x_{c}^2/2}}4.
\end{align*}
It follows then that
\begin{align*}
\frac{I_{22}^{(1)}(x)}{2(1+c)} &:=  (1+x_{c}^2)A_1 - A_2 - A_3 = \frac{1+ x_{c}^2}{\cL{(1+c)^2}}\left[ e^{-x_{c}^2/2}-2\Psi(x_{\cL{c}})\right] \\
&\quad -\frac{[1+x_{c}^2/2]e^{-x_{c}^2/2} - 2[\Psi(x_{\cL{c}})+ x_{\cL{c}}\varphi(x_{\cL{c}})]}{\cL{(1+c)^2}} - \frac {x }4 e^{-x_{c}^2/2} \\
&= \frac{2}{\cL{(1+c)^2}}\left[x_{\cL{c}}\varphi(x_{\cL{c}}) - x_{c}^2 \Psi(x_{\cL{c}})\right].
\end{align*}
\underline{(ii) For $A_4 \sim A_6$ with $c\neq1$}.
We rewrite first $A_4$  and $A_5$ as follows.
\begin{align*}
\cL{c} A_4 &= \int_0^x\Psi(\sqrt{x_{c}^2-t_{\cL{c}}^2})\Psi(|t_{c}'|)e^{-\cL{c}t} \td \cL{c} t\\
&=  \int_0^{cx}\int_{x_{\cL{c}}\sqrt{1-t/(cx)}}^\IF\int_{|x_{\cL{c}}'|\sqrt{t/(cx)}}^\IF e^{-t}\varphi(u)\varphi(v) \td t \td u \td v\\
&=  x_{\cL{c}} |x_{\cL{c}}'|\int_0^\IF\int_0^\IF \varphi(x_{\cL{c}} u)\varphi(|x_{\cL{c}}'|v) \int_{\cL{c}x\max(0, 1-u^2)}^{\cL{c}x\min(1, v^2)} \td [-e^{-t}] \td u \td v, \\
c^2 x A_5 &= x_c^2 \int_0^x (ct)\Psi(\sqrt{x_{c}^2-t_{\cL{c}}^2})\Psi(|t_{c}'|)e^{-\cL{c}t} \td (ct) \\
&= x_{\cL{c}}^3 |x_{\cL{c}}'|\int_0^\IF\int_0^\IF \varphi(x_{\cL{c}}u)\varphi(|x_{\cL{c}}'|v) \int_{\cL{c}x\max(0, 1-u^2)}^{\cL{c}x\min(1, v^2)} \td [-(1+t)e^{-t}] \td u \td v \\
&=c x_c^2 A_4 + x_{\cL{c}}^3 |x_{\cL{c}}'|\int_0^\IF\int_0^\IF \varphi(x_{\cL{c}}u)\varphi(|x_{\cL{c}}'|v) \int_{\cL{c}x\max(0, 1-u^2)}^{\cL{c}x\min(1, v^2)}  \td [-te^{-t}] \td u \td v \\
&=:c x_c^2 A_4 + B_4.
\end{align*}
We decompose further the integrals $A_4$ and $B_4$ by specifying the integral limits as below
\begin{align*}
\cL{c} A_{4}&= x_{\cL{c}} |x_{\cL{c}}'|\left(\int_0^{1} \int_0^{1} \int_{\cL{c}x (1-u^2)}^{\cL{c}x v^2} + \int_0^{1} \int_{1}^\IF \int_{\cL{c}x (1-u^2)}^{\cL{c}x}\right.\\
&\quad +\left. \int_{1}^\IF\int_0^{1}\int_{0}^{\cL{c}xv^2} + \int_{1}^\IF\int_1^\IF\int_{0}^{\cL{c}x} \right)\varphi(x_{\cL{c}} u)\varphi(|x_{\cL{c}}'|v) \td [-e^{-t}] \td u \td v\\
&=: I_1 + I_2 + I_3 + I_4,\\
B_{4}&= x_{\cL{c}}^3 |x_{\cL{c}}'|\left(\int_0^1 \int_0^1 \int_{\cL{c}x(1-u^2)}^{ \cL{c}x v^2 } + \int_0^1 \int_1^\IF \int_{\cL{c}x(1-u^2)}^{\cL{c}x}\right.\\
&\quad +\left. \int_{1}^\IF\int_0^1\int_{0}^{\cL{c}xv^2} + \int_{1}^\IF\int_1^\IF\int_{0}^{\cL{c}x} \right) \varphi(x_{\cL{c}} u)\varphi(|x_{\cL{c}}'|v)\td [-te^{-t}] \td u \td v \\
&=: J_1 + J_2 + J_3 + J_4.
\end{align*}
We deal first with $I_3, I_4$ and $J_3, J_4$.  Since $\varphi(|x_{\cL{c}}|u) e^{cxu^2} = \varphi(|x_{\cL{c}}'|u), \varphi(|x_{\cL{c}}'|v) e^{cxv^2} = \varphi(|x_{\cL{c}}|v)$, we have
\begin{align*}
I_3 
&={\Psi(x_{\cL{c}})}\int_0^1 |x_{\cL{c}}'| \left[\varphi(|x_{\cL{c}}'|v) - \varphi(x_{\cL{c}} v) \right]\td v\\
&= {\Psi(x_{\cL{c}})} \left[ [0.5-\Psi(|x_{\cL{c}}'|)]- \frac {\cL{|1-c|}}{1+c}[0.5-\Psi(x_{\cL{c}})]\right] \\
J_3  
&=-\cL{c}x{\Psi(x_{\cL{c}})}x_{c}^2|x_{\cL{c}}'|\int_0^1 v^2 \varphi(x_{\cL{c}}v) \td v =- \frac{\cL{|1-c|}}{\cL{1+c}} (\cL{c}x){\Psi(x_{\cL{c}})} \omega(x_{\cL{c}})
\end{align*}
and
\begin{align*}
I_4  = 
\Psi(x_{\cL{c}})\Psi(|x_{\cL{c}}'|)[1-e^{-\cL{c}x}], \quad J_4 = 
 -(\cL{c}x )e^{-\cL{c}x} x_{c}^2\Psi(x_{\cL{c}})\Psi(|x_{\cL{c}}'|).
\end{align*}
For $I_1, I_2$ and $J_1, J_2$, we have
\begin{align*}
I_1 
&= x_{\cL{c}}|x_{\cL{c}}'| \int_{0<u, v< 1,u^2+v^2 >1}\left[e^{-\cL{c}x}\varphi(|x_{\cL{c}}'|u)\varphi(|x_{\cL{c}}'|v) - \varphi(x_{\cL{c}}u)\varphi(x_{\cL{c}}v)\right] \td u \td v\\
&= \frac{\cL{1+c}}{\cL{|1-c|}} e^{-\cL{c}x} \mu(|x_{\cL{c}}'|) - \frac{\cL{|1-c|}}{\cL{1+c}} \mu(x_{\cL{c}}), \\
J_1 
&= \cL{c}xx_{\cL{c}}^3 |x_{\cL{c}}'| \int_{0<u, v< 1,u^2+v^2 >1}  \Big[e^{-\cL{c}x}(1-u^2)\varphi(|x_{\cL{c}}'|u)\varphi(|x_{\cL{c}}'|v) - v^2\varphi(x_{\cL{c}}u)\varphi(x_{\cL{c}}v)\Big] \td u \td v\\
&= \cL{c}x \left\{\frac{1+c}{\cL{|1-c|}}x_{c}^2 e^{-\cL{c}x}\mu(|x_{\cL{c}}'|) - \frac{\cL{(1+c)^3}}{\cL{|1-c|}^3} e^{-\cL{c}x} \nu(|x_{\cL{c}}'|) - \frac{\cL{|1-c|}}{1+c} \nu(x_{\cL{c}})\right\},
\end{align*}
where $\mu(a)$ and $\nu(a)$ are given by
\begin{eqnarray}
\begin{split}\label{Omega}
\mu(a)&:= a^2\int_{0<u, v<1,u^2 + v^2 >1} \varphi(a u)\varphi(a v) \td u \td v
\COM{\notag \\
&=  a^2\left(\int_{0<u,v<1} - \int_{u^2+v^2<1, u,v>0}\right)\varphi(a u)\varphi(a v)\td u \td v,\notag \\
&= [\Phi(a) - 0.5]^2 -\frac14[1-e^{-a^2/2}] = \frac{e^{-a^2/2}}{4} - \Phi(a) \Psi(a), \notag\\
\omega(a) &= a^3 \int_0^1 u^2\varphi(au) \td u
} = \Phi(a)-0.5 - a\varphi(a), \\
\nu(a) &:= a^4\int_{0<u, v<1,u^2 + v^2 >1} u^2\varphi(a u)\varphi(a v) \td u \td v
= [\Phi(a)-0.5]\omega(a) +\frac{(1+a^2/2)e^{-a^2/2}-1}4.\notag
\end{split}
\end{eqnarray}
Now, for $I_2$ and $J_2$, we have
\begin{align*}
I_2 
&=\Psi(|x_{\cL{c}}'|) e^{-\cL{c}x} \int_{0}^{1}x_{\cL{c}} \left[\varphi(|x_{\cL{c}}'| u) - \varphi(x_{\cL{c}}u)\right] \td u\\
&=\Psi(|x_{\cL{c}}'|) e^{-\cL{c}x}\left[\frac{\cL{1+c}}{\cL{|1-c|}}[0.5-\Psi(|x_{\cL{c}}'|)] - [0.5-\Psi(x_{\cL{c}})]\right], \\
J_2 
&= \cL{c}xe^{-\cL{c}x}\Psi(|x_{\cL{c}}'|) \int_0^1 x_{\cL{c}}^3\left[(1-u^2)\varphi(|x_{\cL{c}}'|u)- \varphi(x_{\cL{c}}u)\right]\td u  \\
&= \cL{c}xe^{-\cL{c}x}\Psi(|x_{\cL{c}}'|)\left\{\frac{\cL{1+c}}{\cL{|1-c|}} x_{c}^2 [\Phi(|x_{\cL{c}}'|)-0.5] -\frac{\cL{(1+c)^3}}{\cL{|1-c|}^3}\omega(|x_{\cL{c}}'|) - x_{c}^2[\Phi(x_{\cL{c}}) - 0.5] \right\}.
\end{align*}

Therefore, summing up $I_1\sim I_4$ and $J_1\sim J_4$ to give $A_4$ and $B_4$ by putting  $\mu(\cdot),\nu(\cdot),\omega(\cdot)$ together with the fact that $\nu(a) = -0.5 \mathrm m(a) - \Psi(a)\omega(a)$ for $\mathrm m(\cdot)$ given by \eqref{A2}, we have
\begin{align*}
\cL{c} A_{4} &= I_1 + I_2 + I_3 + I_4 \\
&= \frac{\cL{c}}{\cL{|1-c^2|}}e^{-x_{c}^2/2} + \frac{\cL{|1-c|}}{\cL{1+c}} \Phi(x_{\cL{c}})\Psi(x_{\cL{c}}) - \frac{\cL{(1+c)} e^{-\cL{c}x}}{\cL{|1-c|}} \Phi(|x_{\cL{c}}'|)\Psi(|x_{\cL{c}}'|) \\
&\quad +\Psi(|x_{\cL{c}}'|) e^{-\cL{c}x}\left[\frac{\cL{1+c}}{\cL{|1-c|}}[0.5-\Psi(|x_{\cL{c}}'|)] - [0.5-\Psi(x_{\cL{c}})]\right] \\
&\quad + {\Psi(x_{\cL{c}})} \left[ [0.5-\Psi(|x_{\cL{c}}'|)]- \frac {\cL{|1-c|}}{1+c}[0.5-\Psi(x_{\cL{c}})]\right]  + \Psi(x_{\cL{c}})\Psi(|x_{\cL{c}}'|)[1-e^{-\cL{c}x}] \\
&= \max(1,c)\left[\frac{\Psi(x_{\cL{c}})}{1+c}- \frac{\Psi(|x_{\cL{c}}'|)}{\cL{|1-c|}}e^{-\cL{c}x}\right] + \frac{\cL{c}}{\cL{|1- c^2|}} e^{-x_{c}^2/2},\\
 B_4/[\cL{c}x]  &=  (J_1 + J_2 + J_3 + J_4)/(\cL{c}x) =
 \frac{x_{c}^2}{4}\frac{\cL{1+c} }{\cL{|1-c|}} e^{-x_{c}^2/2}  - x_{c}^2\frac{\max(1,c) }{\cL{|1-c|}} e^{-\cL{c}x}\Psi(|x_{\cL{c}}'|) \\
&\quad +\frac{\cL{|1-c|}}{2\cL{(1+c)}}
\left(\Psi(x_{\cL{c}})+x_{\cL{c}}\varphi(x_{\cL{c}})-\frac{(1+x_{c}^2/2)e^{-x_{c}^2/2}}2 \right) \\
&\quad +\frac{\cL{(1+c)^3}}{2\cL{|1-c|}^3} e^{-\cL{c}x}
\left( \Psi(|x_{\cL{c}}'|)  +|x_{\cL{c}}'|\varphi(|x_{\cL{c}}'|)- \frac{(1+|x_{\cL{c}}'|^2/2)e^{-|x_{\cL{c}}'|^2/2}}2 \right)
\\
&=  \frac{x_{c}^2}{4}\frac{\cL{1+c} }{\cL{|1-c|}} e^{-x_{c}^2/2}  - x_{c}^2\frac{\max(1,c) }{\cL{|1-c|}} e^{-\cL{c}x}\Psi(|x_{\cL{c}}'|) + \frac{\cL{|1-c|}}{2\cL{(1+c)}}\mathrm m(x_c) + \frac{\cL{(1+c)^3}}{2\cL{|1-c|}^3} e^{-\cL{c}x}
\mathrm m(|x_c'|).
\end{align*}

Next, we turn to $A_6$.  Since $\varphi(\sqrt{x_{c}^2-t_{\cL{c}}^2})e^{c(x-t)} = \varphi(t_c')$, we have by the integral of $A_2$  given in \eqref{A2}
\begin{align*}
A_6 
&= \frac{1+c}{\cL{|1-c|}}e^{-\cL{c}x} \int_0^x |t_{c}'|\varphi(t_{c}')\Psi\left(\sqrt{|x_{\cL{c}}'|^2-t_{c}'^2}\right) \td t \notag\\
&= \frac{4\cL{(1+c)}}{\cL{|1-c|}^3}e^{-\cL{c}x}\int_0^{|x_{\cL{c}}'|} v^2\varphi(v)\Psi\left(\sqrt{|x_{\cL{c}}'|^2-v^2}\right) \td v  = -\frac{2\cL{(1+c)}}{\cL{|1-c|}^3}e^{-\cL{c}x}\mathrm m(|x_c'|).
\end{align*}
Therefore, it follows by $A_5 = x_c^2/(cx) A_4 +B_4 /(c^2x)$ that
\begin{align}\label{I22_two}
\frac{I_{22}^{(2)}(x)}{2\abs{1-c}} &:=  (1+x_{c}^2)A_4 - A_5 - A_6 = \left[1+x_{c}^2\left(1-\frac1{\cL{c}x}\right)\right] A_4 - \frac{1}c\frac{B_4}{cx} - A_6\notag\\
&=  \left[1+x_{c}^2\left(1-\frac1{\cL{c}x}\right)\right]\left\{\frac{\max(1,c)}c\left[\frac{\Psi(x_{\cL{c}})}{1+c} - \frac{\Psi(|x_{\cL{c}}'|)}{\cL{|1-c|}}e^{-\cL{c}x}\right] + \frac{1}{\cL{|1- c^2|}} e^{-x_{c}^2/2}\right\}\notag \\
&\quad -
\frac{x_{c}^2}{4c}\frac{\cL{1+c} }{\cL{|1-c|}} e^{-x_{c}^2/2}  + \frac{x_{c}^2}{c}\frac{\max(1,c) }{\cL{|1-c|}} e^{-\cL{c}x}\Psi(|x_{\cL{c}}'|)  - \frac{\cL{|1-c|}}{2c\cL{(1+c)}}\mathrm m(x_c)   -   \frac{\cL{1+c}}{2c\cL{|1-c|}}e^{-\cL{c}x}\mathrm m(|x_c'|)\notag\\
&= \left[1+x_{c}^2\left(1-\frac1{\cL{c}x}\right)\right]\frac{\max(1,c)}c \frac{\Psi(x_{\cL{c}})}{1+c}
- \frac{1}{2c|1-c|}\left[\frac{(1-c)^2}{2(1+c)}x_c^2 + \frac{1+c^2}{1+c}\right]e^{-x_{c}^2/2}\notag  \\
&\quad + \frac{1+c^2}{2c^2}\frac{\max(1,c)}{\abs{1-c}}e^{-\cL{c}x}\Psi(|x_{\cL{c}}'|) - \frac{\cL{|1-c|}}{2c\cL{(1+c)}}\mathrm m(x_c)  -   \frac{\cL{1+c}}{2c\cL{|1-c|}}e^{-\cL{c}x}\mathrm m(|x_c'|).
\end{align}

\underline{(iii) For $A_4^+\sim A_6^+$ with $c>1$}. Following similar arguments of $A_4 \sim A_6$, we have
\begin{align*}
\cL{c} A_4^+ &= x_{\cL{c}}\int_0^\IF \varphi(x_{\cL{c}}u)\int_{\cL{c}x\max(0, 1-u^2)}^{\cL{c}x}\td [-e^{-t}] \td u\\
&= x_{\cL{c}} e^{-\cL{c}x} \int_{0}^{1}  \left[\varphi(x_{\cL{c}}'u) - \varphi(x_{\cL{c}}u) \right]\td u + (1-e^{-\cL{c}x})\Psi(x_{\cL{c}}) \\
&= \frac{e^{-\cL{c}x}}{\cL{\abs{1-c}}} \big[1 -(1+c)\Psi(|x_{c}'|)\big] +\Psi(x_{\cL{c}}),\\
c^2x A_5^+
&= x_{\cL{c}}^3\int_0^\IF \varphi(x_{\cL{c}}u)\int_{\cL{c}x\max(0, 1-u^2)}^{\cL{c}x}\td [-(1+t)e^{-t}] \td u =: cx_c^2 A_4^+ + B_4^+,
\end{align*}
with
\begin{align*}
B_4^+ &= x_{\cL{c}}^3\int_0^\IF \varphi(x_{\cL{c}}u)\int_{\cL{c}x\max(0, 1-u^2)}^{\cL{c}x}\td [-te^{-t}] \td u\\
&= x_{\cL{c}}^3\left(\int_0^1 \int_{\cL{c}x(1-u^2)}^{\cL{c}x} + \int_1^\IF\int_0^{\cL{c}x}\right)\td[-te^{-t}] \varphi(x_{\cL{c}}u)\td u \\
&= \cL{c}xe^{-\cL{c}x}\Big[ \int_0^1 x_{\cL{c}}^3[(1-u^2)\varphi(|x_{\cL{c}}'|u)- \varphi(x_{\cL{c}}u)]\td u - x_{c}^2 \Psi(x_{\cL{c}})\Big]\\
&= \cL{c}xe^{-\cL{c}x}\Big[  \frac{x_{c}^2}{\cL{\abs{1-c}}}-\frac{\cL{1+c}}{\cL{\abs{1-c}}}  x_{c}^2\Psi(|x_{c}'|) - \frac{\cL{(1+c)^3}}{\cL{\abs{1-c}}^3}\omega(|x_{c}'|)\Big].
\end{align*}
For $A_6^+$, by a change of variable $s= t-x$ and $e^{ct}\varphi(t_c) = \varphi(|t_{c}'|)$
\begin{align*}
A_6^+ 
= \frac{1+c}{\cL{\abs{1-c}}} e^{-\cL{c}x}\int_0^x |t_{c}'|\varphi(|t_{c}'|)\td t =\frac{4\cL{(1+c)}}{\cL{\abs{1-c}}^3} e^{-\cL{c}x} \int_0^{|x_{\cL{c}}'|}s^2\varphi(s)\td t = \frac{4\cL{(1+c)}}{\cL{\abs{1-c}}^3} e^{-\cL{c}x} \omega(|x_{c}'|).
\end{align*}
Therefore, we have  (recall that $A_5^+ = x_c^2/(cx) A_4^+ +B_4^+ /(c^2x)$)
\begin{align}\label{I22_three}
\frac{I_{22}^{(3)}(x)}{2\abs{1-c}} &= [(1+x_{c}^2)A_4^+ - A_5^+ - A_6^+] \notag\\
&= \left[1+x_{c}^2\left(1-\frac1{\cL{c}x}\right)\right] A_4^+ - \frac1{c^2x}B_4^+ - A_6^+ \notag\\
&=  \left[1+x_{c}^2\left(1-\frac1{\cL{c}x}\right)\right] \left[ \frac{e^{-\cL{c}x}}{\cL{c} } \left(\frac1{\cL{\abs{1-c}}} -\frac{\cL{1+c}}{\cL{\abs{1-c}}} \Psi(|x_{c}'|)\right) +\frac{\Psi(x_{\cL{c}})}c\right] \notag\\
&\quad - \frac{e^{-\cL{c}x}}{\cL{c} }\Big[  \frac{x_{c}^2}{\cL{\abs{1-c}}}-\frac{\cL{1+c}}{\cL{\abs{1-c}}}  x_{c}^2\Psi(|x_{c}'|) - \frac{\cL{1+c}}{\cL{\abs{1-c}}}\omega(|x_{c}'|)\Big]\notag \\
&= - \frac{1+c^2}{2c|1-c|} \frac{e^{-cx}}c \big[1-(1+c)\Psi(|x_{c}'|)\big] + \left[1+x_{c}^2\left(1-\frac1{\cL{c}x}\right)\right]\frac{\Psi(x_{\cL{c}})}c
+  \frac{e^{-cx}}c\frac{\cL{1+c}}{\cL{\abs{1-c}}}\omega(|x_{c}'|).
\end{align}
Consequently, the desired expression of ${\MB}_{\alpha}^h(x)$ follows by \eqref{Ix}, \eqref{def_II2} and \eqref{II1}. Indeed,
\begin{align}\label{c_one}
c[I_1(x) + I_{21}(x) + I_{22}^{(1)}(x)] &=  2\left[ (1+x_{c}^2)\Psi(x_{\cL{c}})-x_{\cL{c}}\varphi(x_{\cL{c}})\right] \notag\\
&\quad + \left[(1+c)\Psi(x_{\cL{c}}) -(1-c)\Psi(x_{c}')e^{-\cL{c}x}\right]
+  \frac{4c}{\cL{1+c}}\left[x_{\cL{c}}\varphi(x_{\cL{c}}) - x_{c}^2 \Psi(x_{\cL{c}})\right] \notag \\
&=  [3+c +(1-c^2)x]\Psi(x_c) -\frac{2(1-c)}{1+c}x_c\varphi(x_c) - (1-c)\Psi(x_{c}')e^{-\cL{c}x},
\end{align}
which together with \eqref{I22_two} implies that
\begin{align}\label{c_two}
\lefteqn{c[I_1(x) + I_{21}(x) + I_{22}^{(1)}(x)] - c I_{22}^{(2)}(x)}\notag\\
&=   [3+c +(1-c^2)x]\Psi(x_c) -\frac{2(1-c)}{1+c}x_c\varphi(x_c) - (1-c)\Psi(x_{c}')e^{-\cL{c}x} \notag\\
&\quad - \left[1+x_{c}^2\left(1-\frac1{\cL{c}x}\right)\right]\frac{\max(1,c)}c \frac{\Psi(x_{\cL{c}})}{1+c}
+\left[\frac{(1-c)^2}{2(1+c)}x_c^2 + \frac{1+c^2}{1+c}\right]e^{-x_{c}^2/2}\notag \\
&\quad - \frac{1+c^2}{c}\max(1,c) e^{-\cL{c}x}\Psi(|x_{\cL{c}}'|)  +  (1+c)e^{-\cL{c}x}\mathrm m(|x_c'|)+ \frac{\cL{|1-c|^2}}{\cL{1+c}}\mathrm m(x_c)\\
&\stackrel{c<1}{=}  \left[\frac{1+c}{c} +\frac{4c}{1+c}\right]\Psi(x_c) -\frac{2(1-c)}{1+c}x_c\varphi(x_c) -\frac{1+c}c\Psi(x_{c}')e^{-\cL{c}x}\notag   \\
&\quad +\left[\frac{(1-c)^2}{2(1+c)}x_c^2 + \frac{1+c^2}{1+c}\right]e^{-x_{c}^2/2}  +  \frac{\cL{(1-c)^2}}{\cL{1+c}}\mathrm m(x_c)  +  (1+c)e^{-\cL{c}x}\mathrm m(x_c').\label{c_uniform}
\end{align}

Now, for $c>1$, we have by \eqref{c_two}
\begin{align*}
{\MB}_{\alpha}^h(x)  &= c[I_1(x) + I_{21}(x) + I_{22}^{(1)}(x) - I_{22}^{(2)}(x)]  +  c I_{22}^{(3)}(x) \\
&=   [3+c +(1-c^2)x]\Psi(x_c)  -\frac{2(1-c)}{1+c}x_c\varphi(x_c) - (1-c)\Psi(x_{c}')e^{-\cL{c}x} \\
&\quad - \left[1+x_{c}^2\left(1-\frac1{\cL{c}x}\right)\right] \frac{\Psi(x_{\cL{c}})}{1+c}
+\left[\frac{(1-c)^2}{2(1+c)}x_c^2 + \frac{1+c^2}{1+c}\right]e^{-x_{c}^2/2}\\
&\quad - (1+c^2) e^{-\cL{c}x}\Psi(|x_{\cL{c}}'|)  + \frac{\cL{|1-c|^2}}{\cL{1+c}}\mathrm m(x_c) +  (1+c)e^{-\cL{c}x}\mathrm m(|x_c'|) \\
&\quad - \frac{1+c^2}{c} e^{-cx} \big[1-(1+c)\Psi(|x_{c}'|)\big] + 2|1-c|\left[1+x_{c}^2\left(1-\frac1{\cL{c}x}\right)\right]\Psi(x_{\cL{c}}) + 2(1+c)e^{-c x}\omega(|x_{c}'|)  \\
&= \left[\frac{1+c}{c} +\frac{4c}{1+c}\right]\Psi(x_c) -\frac{2(1-c)}{1+c}x_c\varphi(x_c) -\frac{1+c}c \Psi(x_{c}')e^{-\cL{c}x} \\
&\quad +\left[\frac{(1-c)^2}{2(1+c)}x_c^2 + \frac{1+c^2}{1+c}\right]e^{-x_{c}^2/2}
+  \frac{\cL{|1-c|^2}}{\cL{1+c}}\mathrm m(x_c)  +  (1+c)e^{-\cL{c}x}\mathrm m(|x_c'|) + 2(1+c)e^{-cx}\omega(|x_{c}'|),
\end{align*}
which equals the right-hand side of \eqref{c_uniform} since the summand of the last two terms equals $(1+c)e^{-\cL{c}x}\mathrm m(x_c') $ using $\mathrm m(|a|) +2\omega(|a|) = \mathrm m(a), a=x_c'<0$ (recall \eqref{A2} and \eqref{Omega}).
Therefore, combining \eqref{c_one} for $c=1$, we have a uniform expression of ${\MB}_{\alpha}^h(x) $  as in \eqref{c_uniform} for all $c>0$. Using again
\begin{align*}
\mathrm m(a) = \Psi(a)+ a\varphi(a) - \frac{[1+a^2/2]e^{-a^2/2}} 2
\end{align*}
and sorting out all terms related to $a=x_c, x_c'$, the desired claim for $\alpha = 1$ follows.

\underline{\textbf{Proof of (ii) $\alpha=2$}}. Recall that $f_\N(t) =(1+c)t^2 -  \sqrt2\N t  + z,\, c>0, t\inr$.  Without loss of generality, we consider only the minimum $f_\N(t^*) <0$ with $t^* =\N/[\sqrt2(1+c)]$. Therefore, there are two solutions   $t_1 < t_2$ of the equation $f_\N(t)=0$  satisfying
\begin{align*}
t_1 = \frac{\sqrt2\N - \sqrt{2\N^2 - 4(1+c)z}}{2(1+c)},\quad t_2 = \frac{\sqrt2\N + \sqrt{2\N^2 - 4(1+c)z}}{2(1+c)}.
\end{align*}
Consequently,
\begin{align*}
{\MB}_2^h(x) &= \int_0^\IF e^z \pk{ \frac{\N}{\sqrt{2}(1+c)}>\frac x2, {t_2-t_1 }>x   } \td z
 + \int_{-\IF}^0 e^z \pk{\frac{\N}{\sqrt{2}(1+c)}< \frac x2, f_\N(x)<0 } \td z\\
&= \int_0^\IF e^z\pk{ 2\N^2 >4(1+c)z+(1+c)^2x^2, \N>0} \td z\\
&\quad +\int_{-\IF}^0 e^z \pk{(1+c)x+z/x<\sqrt2\N <(1+c) x} \td z
=: A(x) + B(x),
\end{align*}
where  it follows by elementary calculations that $B(x) = \Psi\left( (1-c)x/\sqrt{2}\right)e^{-c x^2}$ and
\begin{align*}
A(x) 
= \sqrt{\frac{1+c}{c}}\Psi\left( \sqrt {\frac{c(1+c)}2}x \right)e^{-\frac{(1+c)x^2}{4}} - \Psi\left( \frac{1+c}{\sqrt{2}} x\right).
\end{align*}
The desired claim for $\alpha = 2$ follows.
Consequently, we complete the proof of \netheo{T_one_two}. \QED

Below, we shall verify \netheo{T_two} by noting that $f_\N(t) = [h(t) + t^2] - \sqrt2\N t +z$ for given $z\inr$ is continuous  and convex if and only if $h(t) + t^2$ is. Recall $\mathcal L_{\N}(z)= \int_E \I{f_\N(t)<0}\td t$.

\prooftheo{T_two}
(i) We decompose $\{\mathcal L_{\N}(z) >x\}$ according to the signs of $f_\N(a)$  and $f_\N(b)$. Given $z\inr$, since  $f_\N({s}), {s}\in E$ is continuous and convex, we have  $f_\N(a)<0, f_\N(b)<0$ implies that
$f_\N({s}) <0, {s}\in [a,b]$. Similar arguments for the  cases with $f_\N(a)f_\N(b)<0$. Finally, for $f_\N(a), f_\N(b)>0$, suppose without loss of generality that the minimum of $f_\N(s), s\in E$  is negative. Therefore, it follows by the convexity of $f_\N({s}), {s}\in E$ that, there \ree{exist} two different roots $s_1<s_2$ of $f_\N(s) = 0, s\in E$. Therefore, $\{s\in E: f_\N(s) <0\} = (s_1, s_2)$. Consequently, the first desired claim follows.

(ii) It follows by the convexity of $f_\N(t)$  that $f_\N'(a) := h'_+(a) + 2a - \sqrt2\N \le f_\N'(b) := h'_-(b) + 2b - \sqrt2\N$.  Therefore, we decompose $\{\mathcal L_{\N}(z) >x\}$ according to the three cases that $0\le f_\N'(a)< f_\N'(b)$, $f_\N'(a)<0< f_\N'(b)$ and  $f_\N'(a)< f_\N'(b) \le0$.
For the first case with $f_\N'(a) \ge0$,  it follows by the convexity of $f_\N(t), t\inr$ that
\begin{align*}
\frac{f_\N(t) - f_\N(a)}{t-a}\ge f_\N'(a) \ge0, \quad t>a,
\end{align*}
implying that $\{\mathcal L_\N(z) >x, f_\N'(a) >0\} =\{f_\N(a+x) <0, f_\N'(a) >0\}$. Similar arguments apply for the other two cases.  We complete the proof of \netheo{T_two}.
\QED

\subsection{Proofs of Propositions \ref{T_12} and \ref{T_34}}
\ \newline
In \neprop{T_12}, we take $h(t) = c |t|^\lambda -t^2$ with $\lambda=1,2$}. The general case with $\lambda\ge1$ is shown in \neprop{T_34}. \\
\underline{(i) For $\lambda=2$}. Suppose without loss of generality that the minimum of $f_\N(t), t\inr$ is negative and the two solutions $t_1 < t_2$ of  $f_\N(t)=0, {t}\inr$  satisfying
\begin{align*}
t_1 = \frac{\sqrt2\N - \sqrt{2\N^2 - 4c z}}{2c},\quad t_2 = \frac{\sqrt2\N + \sqrt{2\N^2 - 4c z}}{2c}.
\end{align*}
We have thus  (recall \eqref{MB_E})
\begin{align*}
\mathcal L_\N(z)
=  \min(b, s_2) - \max(a, s_1)
 = \min({b-a, t_2 -a,  b-t_1, t_2-t_1})
\end{align*}
implying that ({note that} $b-a\rr{\ge} x\geq0$)
\begin{align*}
 \MB_2^h(x,{E})
= \int_\R e^z\pk{t_2-a>x, b-t_1>x, t_2-t_1
>x}   \td z = \int_\R e^z \pk{\sqrt{2\N^2-4c z}> u}  \td z,
\end{align*}
with
\begin{align*}
u&= \max( 2c(a+x)-\sqrt{2}\N, \sqrt{2}\N-2c(b-x), c x)\\
&= \left\{
\begin{array}{ll}
2c(a+x)-\sqrt{2}\N, & {\sqrt2\N /c < 2 a + x},\\
\sqrt{2}\N-2c(b-x), &  {\sqrt2\N /c > 2 b - x},\\
c x, & \mbox{otherwise}.
\end{array}
\right.
\end{align*}
Therefore, the remaining argument follows by  elementary calculations and the claim is obtained.

\underline{{(ii) For $\lambda=1$}}.
We have $f_\N({t}) = c|{t}| - \sqrt2\N {t} + z$, a  continuous and piece-wise linear function such that
\begin{align*}
f_\N(0) = z, \quad f_\N'({t}) = c \mathrm{sign}({t}) - \sqrt 2\N, \quad {{t}\neq0}.
\end{align*}
Therefore, we consider below the two cases with $0\le a < b$ and  $a<0<b$.

{As $0\le a < b$}. Clearly, the function $f_\N({t}), {t}\ge0$ {is linear with slope $f_\N'(a) =  f_\N'(b) = c-\sqrt2\N$}. Recalling $b-a\rr{\ge}x$ and {$f_\N(0) = z$}, we have with $\nu(\cdot, \cdot) $ given in the theorem
\begin{align*}
 \MB_{2}^h(x,{E}) &=  \int_\R e^z\pk{f_\N(a+x)<0, f_\N'(a)>0} \td z + \int_\R e^z\pk{ f_\N(b-x)<0, f_\N'(b)<0} \td z \\
 &=   \int_{-\IF}^0 e^z\pk{c+\frac z{a+x} <\sqrt2\N<c} \td z +  \int_0^\IF e^z\pk{\sqrt2\N > c+ \frac z{b-x}} \td z + \Psi(c/\sqrt2) \\
&=  \int_{-\IF}^0 e^z\pk{\sqrt2\N > c+\frac z{a+x}} \td z +  \int_0^\IF e^z\pk{\sqrt2\N > c+ \frac z{b-x}} \td z\\
&= \nu(b-x, c) - \nu(a+x,c) +e^{(a+x)(a+x-c)}
\end{align*}
following elementary calculations.

{As  $a<0<b$}. It follows by \netheo{T_two} (ii) that, with slope {$-c -
\sqrt2\N = f_\N'(a) < f_\N'(b) =c-\sqrt2\N$}
\begin{align*}
\{\mathcal L_\N(z) > x\}  &=  \{ f_\N(a+x)<0, f_\N'(a)>0\} + \{ f_\N(b-x)<0, f_\N'(b)<0\}\\
 &\quad + \{f_\N'(a)<0<  f_\N'(b),  \min(b, t_2) - \max(a, t_1) >x, z<0\},
\end{align*}
where $t_1<0<t_2$ are given by \neprop{T_12}. For the purpose of the explicit expressions of the events involved, we rewrite the first integral based on the  sign of $a+x$
	\begin{align*}
	H_1&= \int_\R e^z\pk{f_\N(a+x)<0, f_\N'(a)>0} \td z \\
	&= \I{a+x<0} \int_\R e^z\pk{\sqrt2\N > \max\left(c, c-\frac z{a+x}\right)} \td z
	+ \I{a+x=0} \Phi(c/\sqrt2) \\
	&\quad + \I{a+x >0} \int_{-\IF}^{-2c(a+x)} e^z\pk{{c}+\frac{z}{a+x}<\sqrt2\N <-c} \td z =  \nu(-(a+x),c).
	\end{align*}
	Similarly, the second integral satisfies
$H_2 = \int_\R e^z\pk{f_\N(b-x)<0, f_\N'(b)<0} \td z = \nu(b-x,c).$
Consequently, the desired claims of \neprop{T_12} are obtained.
\QED

In view of \eqref{MB_E}, we have $f_\N(s) = c|t|^\lambda -\sqrt2\N t + z$ is convex such that
\begin{align}
\label{derivative_2}
\qquad\ f_\N(0) = z, \quad  f'_\N(t) = \lambda c  \abs{t}^{\lambda - 1} \mathrm{sign}(t)- \sqrt 2 \N,\quad t\neq0.
\end{align}

\proofprop{T_34} Clearly, it follows from  \netheo{T_two} (ii) that, the claim for $a<0<b$ holds since  {$h'_+(a) +2a= -\lambda c (-a)^{\lambda-1}<0< h_-'(b) + 2b= \lambda c b^{\lambda-1}$} \ree{exist} with finite values.  For $0\le a<b$,  it follows from \eqref{derivative_2} that
\begin{align*}
\small {\pk{\mathcal L_{\N}(z) >x}} = \pk{f_\N(a+x) < 0,\, z<0} +  \pk{{\min(b, t_2) - \max(a, t_1) >x, f_\N(t^*) <0}, \N >0, z>0}.
\end{align*}
Consequently,
\begin{align*}
\lefteqn{\MB_{2}^h(x,{E}) =
  \int_{-\IF}^{0}e^{z}\pk{f_\N(a+x) < 0} \td z }   \\
&\quad +\int_0^{\IF} e^z\pk{{\min(b, t_2) -  \max(a, t_1) >x,\, f_\N(t^*) <0}, \N >0}   \td z,
\end{align*}
where the first integral is $\E{\mathbb I(\sqrt 2\N (a+x) + \W> c(a+x)^{\lambda})}$ by a change of variable $z' = -z$. We compete the proof of \neprop{T_34}.
\QED

\subsection{Proof of Proposition \ref{T5} for the bounds of $\MB_\alpha^h$}
\ \newline
We start with the drift function $h(t) = c|t|^\lambda - t^2$ and $\alpha=2$, and then for the general bounds available for $\alpha\in(0,2]$.\\
 \underline{(i) Bounds of $\MB_2^h$ for  $h(t) = c|t|^\lambda - t^2$ with $\lambda\ge1$}. Clearly, $g(t)=c\abs{t}^\lambda, t\in[a,b]$ and $\lambda\ge1$ is convex and thus
\begin{align}\label{Ineq_g1}
g(t) \ge K_y(t) := g(y) + g'(y)(t-y),\quad  t\in {E}\setminus\{0\}.
\end{align} Therefore, with $c_0 := c_0(y) = g'(y){=\lambda c \abs{y}^{\lambda-1}\sign(y)}, y\in {E}\setminus\{0\}$
\begin{align*}
\MB_2^h(x,{E}) &\le   \int_\R e^z \pk{\int_a^{b} \mathbb I( \sqrt2\N t- K_y(t) -z > 0) \td t > x} \td z \notag \\
&=  \int_\R e^z \pk{\int_a^{b} \mathbb I( \sqrt2\N t- c_0t -[z-c_0 y+g(y)]>0) \td t > x} \td z \\
 &= \expon{c_0(y) y -g(y)} \int_\R e^z \pk{\int_a^{b}\mathbb I(c_0 t - \sqrt2\N t + z <0) \td t > x} \td z
 \\
&=: C_0(y) D_0(y, {E}), \quad y\in{E}\setminus\{0\}.
\end{align*}

Next, we verify $D_0(y)$ satisfies \eqref{Eq_C0} by considering simply  $\tilde f_\N(t) = (c_0 - \sqrt2\N)t + z$. Indeed, we have $\tilde f_\N'(a) = c_0 - \sqrt2\N, c_0\inr$ and $\tilde f_\N(0) = z$. Therefore,
\begin{align*}
D_0(y, {E}) = \int_{\R} e^z \pk{\tilde f_\N(a+x) <0, \tilde f_\N'(a)>0} \td z + \int_{\R} e^z \pk{\tilde f_\N(b-x) <0, \tilde f_\N'(a)<0} \td z,
\end{align*}
where the first integral equals
\begin{align*}
 &\quad \I{a+x{\le}0}\Phi(c_0/\sqrt2) + \I{a+x<0} \int_{0}^\IF e^z\pk{\sqrt2\N <c_0+\frac{z}{a+x}} \td z \notag\\
&\quad +\I{a+x>0} \int_{-\IF}^0 e^z\pk{c_0+\frac{z}{a+x}<\sqrt2\N <c_0} \td z \notag\\
&= e^{(a+x)(a+x-c_0)} \Phi\left(\frac{c_0 -2(a+x)}{\sqrt2}\right) = \nu'(-(a+x), -c_0).
\end{align*}
Similarly, the second integral equals $\nu'(b-x, c_0)$. The  claim follows by the arbitrary of $y\in[a, b]\setminus\{0\}$.

\underline{(ii) Bounds for $\MB_\alpha^h(x, E)$ with $\alpha\in(0,2]$}. {Recalling} $M = \max_{s\in E}h(t)$, we have
\begin{align*}
\mathcal L_\alpha^h(z,E) &:=  \int_E\mathbb I(\sqrt2 B_\alpha(t) - |t|^\alpha - h(t) - z>0) \td t\\
&\ge \int_E\mathbb I(\sqrt2 B_\alpha(t) - |t|^\alpha - [z + M] >0) \td t = \mathcal L_\alpha(z+M,E).
\end{align*}
Hence,
\begin{align*}
\MB_\alpha^h(x, E)=\int_\R e^z \pk{\mathcal L_\alpha^h(z, E) >x}  \td z
\ge \int_\R e^z \pk{\mathcal L_\alpha(z+M, E) >x}  \td z  = e^{-M} \MB_\alpha(x, E).
\end{align*}
Conversely,  the sojourn time $\mathcal L_\alpha^h(z,E)$ is increasing with respect to $T$ involved in the time interval $E=[0,T]$.
We have
\begin{align*}
\mathcal L_\alpha^h(z,E) \le  \int_{0}^\IF \mathbb I(\sqrt2 B_\alpha(t) - |t|^\alpha - h(t) - z>0) \td t = \mathcal L_\alpha^h(z,[0,\IF]),
\end{align*}
implying that $\MB_\alpha^h(x, E) \le \MB_\alpha^h(x, [0,\IF])$. Meanwhile, it is clear that $\MB_\alpha^h(x, E)$ is decreasing with respect to $x\ge0$. The claim follows. We complete the proof of \neprop{T5}.
\QED

\vspace{0.5mm}
\noindent {\bf Acknowledgements}~~
C. Ling would like to thank Prof. Krzysztof D\c{e}bicki  for several useful discussions and important comments during the work on the contribution.
C. Ling is supported by the National Natural Science Foundation (NSNF) (11604375). H. Zhang is partially supported by  the NSNF (11701469) and the Basic and Frontier Research Program of Chongqing, China \rr{(cstc2016jcyjA0510)}.

\section{Appendix}\label{Appendix}
In Section \ref{append_1}, we discuss first  $\MB_\alpha^h(x, E)$ with drift function such that $h(t)+t^2$ is continuous and concave in \netheo{Cor1}, which is illustrated by \neprop{Cor2} with $h(t)=c\abs{t}^\lambda-t^2, 0<\lambda<1$ as well as  its lower bounds in \neprop{Cor3}. Second, we present in Section \ref{append_2} for the detailed calculations of \neprop{T_12}.

\subsection{Discussions on $\MB_\alpha^h$ with $h(t)+t^2$ being concave}\label{append_1}
\ \newline
Recall the curve family $f_\N(t) = f_\N(t,z), z\inr$ (recall \eqref{MB_E}) and the sojourn time $\mathcal L_\N(z)$ given by
\begin{align*}
f_\N(t) = h(t) + t^2 - \sqrt 2 \N t + z, \quad t\inr,\quad \mathcal L_{\N}(z) = \int_E \mathbb I(f_\N(t) <0) \td t.
\end{align*}
Thus $f_\N(t)$ is continuous and concave if and only if $h(t) + t^2$ is.  Let thus $s_1 < s_2$ and $t_1 < t_2$  be the two random solutions of $f_\N(s)=0,\, s\in E$ and $f_\N(t)=0, t\inr$, respectively if it holds that
\begin{align*}
f_\N(s^*) = \max_{s\in E}f_\N(s) >0, \quad f_\N(t^*) = \max_{t\inr}f_\N(t) >0.
\end{align*}

\BT\label{Cor1} Let $\MB_\alpha^h(x, E)$ be given by \eqref{MB_E} with $E = [a,b], a<b, a,b\inr$.

(i) If $h(t)+ t^2$ is a continuous, concave function on $E$, then
\begin{align*}
\lefteqn{\MB_2^h(x, E)  =   \int_{\R} e^z\pk{f_\N(s^*) {\le}
	0} \td z}  \\
&\quad + \int_{\R} e^z \big[\pk{f_\N(a+x) <0, f_\N(b)\ge0} + \pk{f_\N(a)\ge0,  f_\N(b-x) <0}\big] \td z\\
&\quad +\int_{\R} e^z \pk{b-a -(s_2 - s_1) >x, f_\N(s^*)>0, f_\N(a)<0, f_\N(b)<0} \td z.
\end{align*}
(ii) If $h(t)+t^2, t\inr$ is continuous and concave, and the finite right derivative $h'_+(a)$ and the left derivative $h'_-(b)$ exist with finite values, then
{$h'_+(b) \le h'_-(a)$} and
\begin{align*}
\lefteqn{
	\MB_2^h(x, E) = \int_{-\IF}^{-h(0)} e^z\pk{f_\N(t^*)  {\le}
		0} \td z
} \\
&\quad +\int_\R e^z \left[ \pk{f_\N(a+x)<0, \sqrt2\N \le  h'_-(b) + 2b} + \pk{f_\N(b-x)<0, \sqrt2\N \ge h'_+(a) + 2a} \right] \td z \\
&\quad +\int_{\R} e^z\mathbb{P}\{\max({t_1}-a,0) + \max(b-{t_2},0) >x,  h'_-(b) + 2b < \sqrt2\N < h'_+(a) + 2a, f_\N({t^*})>0\} \td z.
\end{align*}
\ET

\textit{Proof}.
(i) We decompose $\{\mathcal L_{\N}(z) >x\}$ according to the signs of $f_\N(a)$  and $f_\N(b)$. Given $z\inr$, since  $f_\N({s}), {s}\in E$ is continuous and concave, we have the maximum $f_\N(s^*) \le 0, s^* \in E$ implies that
$f_\N({s}) \le 0, {s}\in E$. While $f_\N(s^*) > 0$, we analyze the three cases with $f_\N(a) <0, f_\N(b)\ge0$;  $f_\N(a) \ge0, f_\N(b)<0$; and $f_\N(a), f_\N(b) <0$.  For the first case, we have by the concavity
\begin{align*}
\frac{f_\N(t) - f_\N(a)}{t-a} \ge \frac{ f_\N(b) - f_\N(a) }{b-a}\ge0, \quad a < t \le b
\end{align*}
and if there is a $t_0\in[a+x,b)$ such that $f_\N(t_0)=0$, then it holds that $f_\N(t)\ge0$ for all $t\in[a+x,b)$. Therefore, $\{\mathcal L_\N(z) > x, f_\N(a) <0, f_\N(b)\ge0\} = \{f_\N(a+x)<0, f_\N(b)\ge0\}$.
Similar argument applies for the second case. Finally, for $f_\N(a), f_\N(b)<0$ and $f_\N(s^*)>0, s\in E$, it follows by the concavity of $f_\N({s}), {s}\in E$ that, there \ree{exist} two different roots $s_1<s_2$ of $f_\N(s) = 0, s\in E$. Therefore, $\{s\in E: f_\N(s) <0\} = (a, s_1) \cup (s_2, b)$. Consequently, the first desired claim follows.

(ii) For the first case, $f_\N(t^*) $ is non-negative implies that $f_\N({t}) \le 0, t\inr, {z\in (-\IF, -h(0) )}$. The rest cases follow by the concavity of $f_\N(t)$  that $f_\N'(a) := h'_+(a) + 2a - \sqrt2\N \ge f_\N'(b) := h'_-(b) + 2b - \sqrt2\N$.  Therefore, we decompose $\{\mathcal L_{\N}(z) >x\}$ according to the three cases that $f_\N'(b)\ge0$, $f_\N'(b)<0< f_\N'(a)$ and  $f_\N'(a) \le0$.
For the case with $f_\N'(b) \ge0$,  it follows by the concavity of $f_\N(t), t\inr$ that
\begin{align*}
{\frac{f_\N(t) - f_\N(a+x)}{t-(a+x)}\ge f_\N'(b) \ge0, \quad a\le t<a+x,}
\end{align*}
implying that $\{\mathcal L_\N(z) >x, f_\N'(b) \ge 0\} =\{f_\N(a+x) <0, f_\N'(b) \ge 0\}$. Similar arguments apply for the other two cases.  We complete the proof of \netheo{Cor1}.
\QED

\BP\label{Cor2}
Let $\MB_2^h(x,E)$ be  the Piterbarg-Berman function  defined in \eqref{MB_E} with drift function $h(t)= c|t|^\lambda - t^2, c>0, 0 < \lambda < 1$ and $E=[a, b]$.

(i) For $a\ge0$, we have
\begin{align*}
\lefteqn{ \MB_{2}^h(x,E)
= \int_0^{\IF} e^z\pk{f_\N(b-x) < 0} \td z }\\
&\quad + {\int_{-\IF}^0} e^z\pk{ f_\N(a+x) <0, \N<0}   \td z + {\int_{-\IF}^0} e^z\pk{ f_\N(t^*) <0, \N>0}   \td z\\
&\quad + {\int_{-\IF}^0}  e^z\pk{ \max(b-t_2,0)+\max(t_1-a,0)   >x, f_\N(t^*) >0, \N>0 }   \td z.
\end{align*}
(ii) For $a<0$, we have
\begin{align*}
\lefteqn{\MB_{2}^h(x,E)
= \int_{-\IF}^{{-h(0)}} e^z \pk{ f_\N(t^*) <0} \td z} \\
&\quad + \int_\R e^z\left[\pk{f_\N(a+x)<0, f_\N'(b) \ge 0}  + \pk{ f_\N(b-x)<0, f_\N'(a) \le 0} \right] \td z\\
&\quad + \int_\R e^z\pk{\max(t_1-a,0) - \max(b-t_2,0) >x,  f_\N'(b) <0 < f_\N'(a), f_\N(t^*)>0} \td z.
\end{align*}
\cL{Here $t_1< t_2$  are the solutions of $f_\N(t)=0, t\inr$ as its maximum $f_\N(t^*)$ is greater than zero.}
\EP

\textit{Proof}.  Clearly, it follows from  \netheo{Cor1} (ii) that, the claim for $a<0<b$ holds since { $h'_+(a) + 2a = -\lambda c (-a)^{\lambda-1}<0< h_-'(b) + 2b= \lambda c b^{\lambda-1}$ }\ree{exist} with finite values.  For $0\le a<b$, since
\begin{align}
\label{derivative_2}
\qquad\ f_\N(0) = z, \quad  f'_\N(t) = \lambda c  t^{\lambda - 1}- \sqrt 2 \N,\quad t\ge0,
\end{align}
we have
\begin{align*}
 \lefteqn{\pk{\mathcal L_{\N}(z) >x} = \pk{f_\N(b-x)<0, z>0}} \\
&\quad +\pk{f_\N(a+x) < 0,\N<0, z<0} +\pk{ f_\N(t^*) <0, \N>0, z<0} \\
&\quad +\pk{ \max(b-t_2,0)+\max(t_1-a,0)   >x, f_\N(t^*) >0, \N>0, z<0 }.
\end{align*}
 We compete the proof of \neprop{Cor2}.
\QED

\BP
\label{Cor3}
Let $\MB_\alpha^h(x, E)$
be the Piterbarg-Berman function given by \eqref{MB_E}.

 For $\alpha=2$ and $h(t) = c|t|^\lambda - t^2, c>0, 0< \lambda <1$, we have for $E=[a,b]$
\begin{align*}
\MB_2^h(x,E) \ge \left\{ \begin{array}{ll}
	\max_{y\in[a,b]\setminus\{0\}} C_0(y) D_0(y, E), & ab\ge0,\\
	\max\big(\max_{ y\in[a,0)} C_0(y) D_0(y,\zs{[a,0]}),  \max_{ y\in (0, b]} C_0(y) D_0(y,[0,b])\big), & ab<0.
\end{array}
\right.
\end{align*}
Here $C_0(y)$ and $D_0(y)$ are as in \eqref{Eq_C0}. 
\EP

\textit{Proof}.
The main arguments are similar to those for \neprop{T5}. Let $g(t)=c\abs{t}^\lambda, t\in[a,b], 0< \lambda <1$, which is continuous and concave. Therefore, instead of \eqref{Ineq_g1}, we have
\begin{align}\label{Ineq_g2}
g(t) \le K_y(t) := g(y) + g'(y)(t-y),\quad t\in {E},
\end{align}
with {$y\in {E}\setminus\{0\}$ and $yt\ge0$}.

First, for $a\ge0$, we have by a change of variable $z' = z-c_0 y+g(y)$
\begin{align}\label{Same_ab}
\MB_2^h(x,{E}) & \ge  \int_\R e^z \pk{\int_a^{b} \mathbb I( \sqrt2\N t- K_y(t) -z > 0) \td t > x} \td z \notag \\
&=  \int_\R e^z \pk{\int_a^{b} \mathbb I( \sqrt2\N t- c_0t -[z-c_0 y+g(y)]>0) \td t > x} \td z \notag\\
 &= \expon{c_0(y) y -g(y)} \int_\R e^z \pk{\int_a^{b}\mathbb I(c_0 t - \sqrt2\N t + z <0) \td t > x} \td z \notag \\
&=: C_0(y) D_0(y, {E}), \quad y\in{E}\setminus\{0\}.
\end{align}
Next, we deal with the case of  $a<0<b$. Clearly, we have
\begin{align*}
\MB_2^h(x,E) \ge \MB_2^h(x,[a,0])
\ge C_0(y) D_0(y,\zs{[a,0]}),\quad y\in[a,0),
\end{align*}
where the second inequality holds by \eqref{Same_ab}. Similar argument yields that
\begin{align*}
\MB_2^h(x,E)  \ge \MB_2^h(x,[0, b])   \ge C_0(y) D_0(y,\zs{[0,b]}),\quad y\in(0,b].
\end{align*}
Consequently,
\begin{align*}
\MB_2^h(x,E) \ge \max(\max_{ y\in[a,0)} C_0(y) D_0(y, {[a,0])},  \max_{ y\in (0, b]} C_0(y) D_0(y,[0,b]).
\end{align*}

By the arbitrary of $y$, we complete the proof of \neprop{Cor3}.
\QED

\COM{\BP\label{T_12} (Proposition 2.1 in \cite{LZB})
Let {$ \MB_2^h(x,E)$} be  the Piterbarg-Berman function  defined in \eqref{MB_E} with $h(t) = c|t|^\lambda - t^2,c>0$ and $E=[a,b]$. Denote by  $f_\N(t) =  c |t|^\lambda - \sqrt 2 \N t + z,\, \N \sim N(0,1)$ as in \eqref{def_f}.

(i) For $\lambda =2$, we have
\begin{align*}
\small \lefteqn{{\MB_2^h(x,E)}= \small \int_\R e^z \left[\pk{f_\N(a+x) <0, \sqrt 2\N/c< {2a+x}} + \pk{f_\N(b-x)<0, \sqrt 2\N/c>{2b-x}}\right]  \td z}\qquad\qquad\quad  \\
& \small \, + \int_\R e^z \pk{2\N^2/c^2>{ x^2+4 z/c},\,{2a+x}<{\sqrt{2}}\N/c<{2b-x} }  \td z\qquad\qquad\qquad\qquad
\end{align*}
and for $h\equiv0$, i.e., $c=1$
\begin{align*}
\MB_2(x,[0,T]) = 2\Psi(x/\sqrt2)+ \sqrt2(T-x) \varphi(x/\sqrt2). 
\end{align*}
(ii) For $\lambda=1$, we have with $ \nu(m,c) =  e^{m^2 -c|m|} \Psi\left([c-2m]/{\sqrt2}\right)$\begin{align*}
\MB_2^h(x,{E}) = \small \left\{
\begin{array}{ll}
 e^{(a+x)(a+x-c)}  - \nu(a+x,c) + \nu(b-x,c), & a\ge0,\\
 \nu(b-x, c) + \nu(-(a+x), c) +  \int_{-\IF}^0 e^z\pk{\sqrt2|\N|<c,  \min(b, t_2) - \max(a, t_1) >x} \td z,
 & a<0,
  \end{array}
  \right.
\end{align*}
{where $t_1<t_2$ are  the random solutions of $f_\N(t)=0,\, t\inr$ equal}
\begin{align*}
t_1 = \frac z{c+\sqrt2\N}<0< \frac{-z}{c-\sqrt2\N} = t_2,\quad \cL{z<0}.
\end{align*}
\EP
}
\subsection{Explicit expressions of integrals involved in $\MB_2^h (x, {E})$ in \neprop{T_12}} \label{append_2}
\ \newline 
\underline{(i) For $\lambda =2$}. Denote by $f_\N(t) =  c t^2 -  \sqrt 2 \N t + z$ and
  \begin{align*}
I_1(a,x) &= \int_\R e^z \pk{f_\N(a+x) <0, \sqrt 2\N /c< {2a+x}}   \td z, \\
I_2(b,x)&=  \int_\R e^z \pk{f_\N(b-x)<0, \sqrt 2\N/c >{2b-x}}  \td z = I_1(-b,x),\\
I_3(a,b,x) &= \int_\R e^z \pk{2\N^2> {c^2x^2+4c z},\,{2a+x}<{\sqrt{2}}\N /c< {2b-x} }  \td z
\end{align*}
where the relationship between $I_1, I_2$ follows by the symmetry of $\N$.
According to the sign of $a+x$,
\begin{align}\label{eq_I1}
I_1(a,x) &=   \I{a+x>0} \int_{\R} e^z \pk{a+x+\frac z{c(a+x)} <\sqrt 2\N/c< {2a+x}}   \td z \notag \\
&\quad +\I{a+x=0} \pk{\sqrt 2\N/c< {2a+x}}  \notag \\
&\quad +\I{a+x<0} \int_{\R} e^z \pk{\sqrt 2\N/c< \min\left(a+x+\frac z{c(a+x)}, {2a+x}\right)}   \td z.
\COM{ \int_{-\IF}^0 e^z \pk{f_\N(a+x) <0, \sqrt 2\N< c({2a+x}), s_1 <0 <s_2}   \td z \notag\\
&\quad +\I{a+x>0}  \int_{0}^\IF e^z \pk{f_\N(a+x) <0, \sqrt 2\N< c({2a+x}), \N > \sqrt{2z}}   \td z \notag \\
&\quad + \I{a+x<0} \int_{0}^\IF e^z \pk{f_\N(a+x) <0, \sqrt 2\N< c({2a+x}), \N < -\sqrt{2z}}   \td z \notag \\
&=  \I{a+x>0}\int_{-\IF}^{\min(0, a(a+x))} e^z \pk{a+x + \frac{z}{a+x}<\sqrt 2\N< {2a+x}} \td z \notag \\
&\quad +  \I{a+x<0}\int_{-\IF}^{0} e^z \pk{\sqrt 2\N/c< {2a+x}} \td z \notag \\
&\quad +\I{a>0}\int_0^{a(a+x)}e^z \pk{a+x + \frac{z}{a+x}<\sqrt 2\N< {2a+x}} \td z \notag \\
&\quad +\I{a+x<0} \int_{0}^\IF e^z \pk{\sqrt 2\N/c< \min\left(a+x+\frac z{a+x}, {2a+x}\right)}   \td z.
}
\end{align}
Finally, we deal with $I_3$ below. {Rewriting $\{2\N^2/c^2 > x^2+4 z/c\}$ according to $x^2+4z /c\ree{\ge}, <0$}, we have by elementary calculations
\begin{align}\label{eq_I3}
I_3(a,b,x) &= e^{-c x^2/4} \pk{{2a+x}<{\sqrt{2}}\N/c<{2b-x}}\notag \\
&\quad + \int_{- c x^2/4}^\IF e^z \pk{\max(2a+x, \sqrt{x^2+4z/c}) <\sqrt 2\N /c< {2b-x}}   \td z \notag \\
&\quad +  \int_{- c x^2/4}^\IF e^z \pk{2a+x <\sqrt 2\N/c< \min( -\sqrt{x^2+4z/c},{2b-x})}   \td z.
\end{align}
One can further write down the related integrals with respect to $z$ by specifying the maxima/minima involved with restriction of the left-endpoint being less than the corresponding right-endpoint.
Consequently, sum up all related terms given in  \eqref{eq_I1}$\sim$\eqref{eq_I3} and thus the explicit expression of $\MB_2^h(x, {E})$ is obtained.

\CL{\underline{Explicit expression of $\MB_2(x,[0,T])$}.} Replace $a=0, b=T$ and $c=1, x\in[0,T]$ in the expressions of $I_1, I_2$ and $I_3$ above. First,
\begin{align*}
I_1(0, x) =   \I{x>0} \int_{\R} e^z \pk{x+\frac z{x} <\sqrt 2\N < {x}}   \td z  + \I{x=0} \pk{\sqrt 2\N< {x}},
\end{align*}
where the integral equals
\begin{align*}
&\int_{z<0}\int_{x+z/x<\sqrt2 u <x}e^z \varphi(u)\td u   \td z =\int_{\sqrt2u <x}\int_{z<x(\sqrt2 u-x)}e^z \varphi(u)   \td z\td u \\
& = \int_{\sqrt2u <x}e^{x(\sqrt2 u-x)} \varphi(u)\td u =  \int_{\sqrt2(u-\sqrt2x) <-x} \varphi(u - \sqrt2x)\td u = \Psi(x/\sqrt2).
\end{align*}
Consequently, we have $I_1(0, x) =  \Psi(x/\sqrt2), \, x\ge0$.

Now, for $I_2(T,x)$,  we have (recall $-T+x<0$)
\begin{align*}
I_2(T, x) &= I_1(-T, x) = \int_{\R} e^z \pk{\sqrt 2\N< \min\left(-T+x+\frac z{-T+x}, {-2T+x}\right)}   \td z \\
&= \int_{T(T-x)}^\IF e^z \pk{\sqrt 2\N <-T+x+\frac z{-T+x}} \td z +  \int_{-\IF}^{T(T-x)} e^z \pk{\sqrt 2\N <-2T+x} \td z  \\
&= \int_{T(T-x)}^\IF e^z \pk{\sqrt 2\N >T-x+\frac z{T-x}} \td z +  e^{T(T-x)} \Psi((2T-x)/\sqrt2),
\end{align*}
where the first integral is (set $m:=T-x$)
\begin{align*}
 &\quad \int_{\sqrt2u>2m+x} [e^{m(\sqrt2u-m)} - e^{m(m+x)}] \varphi(u) du \\
 &= \int_{\sqrt2u>2m+x} \varphi(u-\sqrt2m) -e^{m(m+x)} \varphi(u) du\\
 &= \Psi(x/\sqrt2) - e^{m(m+x)} \Psi((2m+x)/\sqrt2)).
\end{align*}
Thus, $I_2(T,x) =  \Psi(x/\sqrt2), \, x\ge0$.

Now, we analyze $I_3$ as below
\begin{align*}
I_3(0,T,x) &= e^{-x^2/4} \pk{{x}<{\sqrt{2}}\N<{2T-x}} + \int_{-  x^2/4}^\IF e^z \pk{\max(x, \sqrt{x^2+4z}) <\sqrt 2\N< {2T-x}}   \td z \\
&=\pk{{x}<{\sqrt{2}}\N<{2T-x}}+ \int_0^\IF e^z \pk{\sqrt{x^2+4z} <\sqrt 2\N< {2T-x}}   \td z\\
 &=\pk{{x}<{\sqrt{2}}\N<{2T-x}} + \int_{x<\sqrt2u< 2T-x} \int_{0<4z< 2u^2-x^2} \varphi(u)e^z \td z \td u\\
 &=  \int_{x<\sqrt2u< 2T-x}e^{(2u^2-x^2)/4} \varphi(u)\td u= \sqrt2(T-x) \varphi(x/\sqrt2),
\end{align*}
\CL{which together with $I_1, I_2$ implies the desired result of $\MB_2(x, [0,T])$.
}

\underline{(ii) For $\lambda =1$}. We specify the following integral
\begin{align*}
H_3 =  \int_{-\IF}^0 e^z\pk{{\sqrt 2 |\N| < c},  \min(b, s_2) - \max(a, s_1) >x} \td z. 
\end{align*}
\COM{We put the expressions of $s_1, s_2$ and then rewrite $H_3$ as follows.
\begin{align*}
&\int_{-\IF}^0 e^z\mathbb P\left\{-c < \sqrt2\N <c, (\sqrt2\N -c)(a+x) - z >0, \right. \\
& \quad (\sqrt2\N +c)(b-x) - z >0, 2\N^2 >c^2 + \frac{2cz}x\Big\}\td z \\
&= \pk{\W >\max\left( (c-\sqrt2\N )(a+x), (c+\sqrt2\N)(-b+x), (c^2-2\N^2)\frac x{2c}, c^2>2\N^2\right)}
\end{align*}
}
{First, we rewrite $H_3$ by comparing $t_i,\,i=1,2$ and $a,b$, we have
	\begin{align}\label{eq_H3}
	\zs{H_3} &= \int_{-\IF}^0 e^z \pk{-c<\sqrt2\N<c,  b - t_1>x, a<t_1<b<t_2} \td z \notag\\
	&\quad +  \int_{-\IF}^0 e^z  \pk{-c<\sqrt2\N<c,  t_2 - a>x,t_1 < a <t_2 < b} \td z\notag \\
	&\quad +  \int_{-\IF}^0 e^z  \pk{-c<\sqrt2\N<c, t_1<a< b<t_2} \td z \notag\\
		&\quad + \int_{-\IF}^0 e^z \pk{-c<\sqrt2\N<c,  t_2 - t_1 > x,a<t_1 < t_2 <b} \td z \notag\\
	&=: H_{31} + H_{32} + H_{33} + H_{34},
	\end{align}
	\zs{where
		\begin{align*}t_1 = \frac z{c+\sqrt2\N}<0< t_2 = \frac{-z}{c-\sqrt2\N}, \quad {z<0}.
		\end{align*}}
We deal with the four parts $H_{31}\sim H_{34}$ one by one. Indeed,	
\begin{align*}
\small H_{31} &=  \int_{-\IF}^0 e^z \pk{-c<\sqrt2\N<c, f_\N(a)>0, f_\N(b-x)<0, f_\N(b)<0} \td z \\
&=  \I{b-x{\ge}0}\int_{-\IF}^0 e^z \pk{-c<\sqrt2\N<c, f_\N(a)>0, f_\N(b)<0} \td z \\
&\quad + \I{b-x<0}\int_{-\IF}^0 e^z \pk{-c<\sqrt2\N<c, f_\N(a)>0,  f_\N(b-x)<0, f_\N(b)<0} \td z  \\
\small &=  \small \I{b-x{\ge}0} \int_{-\IF}^0 e^z \pk{\beta<\sqrt2\N<c} \td z \\
&\quad +\I{b-x<0}\int_{-\IF}^0 e^z \pk{\beta<\sqrt2\N<\min\left(c, \frac{z}{b-x}-c\right)} \td z,
\end{align*}
where $\beta$ is given by
$$ \beta  :=\max \left( \frac{z}{a}-c, \frac{z}{b}+c \right)=
\left\{
\begin{array}{ll}
\frac{z}{a}-c, & z<\frac{2abc}{b-a},\\
\frac{z}{b}+c , &  \frac{2abc}{b-a} <z <0.
\end{array}
\right. $$
Similar argument for $H_{32}$ implies that (set below $\alpha = z/a + z/b - \beta$)
\begin{align*}
H_{32} 
&= \I{a+x>0}  \int_{-\IF}^0 e^z \pk{\max\left(-c, \frac{z}{a+x}+c\right)<\sqrt2\N<\alpha} \td z\\
&\quad +\I{a+x{\le}0}  \int_{-\IF}^0 e^z \pk{-c<\sqrt2\N<\alpha} \td z,
\end{align*}
Next, for $H_{33}$, we have
\begin{align*}
H_{33} &=  \int_{-\IF}^0 e^z \pk{-c<\sqrt2\N<c, f_\N(a)<0, f_\N(b)<0} \td z\\
&= \int_{-\IF}^\frac{2abc}{b-a}e^z \pk{-c<\sqrt2\N<c, \frac{z}b+c < \sqrt{2}\N< \frac{z}{a} -c} \td z
\end{align*}
and similar argument for $H_{34}$ implies that
\begin{align*}
H_{34} 
&=  \int_{\frac{2abc}{b-a}}^0 e^z \pk{c^2+\frac{2cz}{x}<2\N^2<c^2,  \frac{z}{a} -c < \sqrt{2}\N<  \frac{z}b+c} \td z\\
&= \I{x< -\frac{4ab}{b-a}} \int_{\frac{2abc}{b-a}}^{-cx/2}\pk{ \frac{z}{a} -c < \sqrt{2}\N<  \frac{z}b+c} \td z \\
&\quad + \int_{\frac{2abc}{b-a}}^0 e^z \pk{\max\left(\sqrt{c^2+\frac{2cz}{x}},  \frac{z}{a} -c\right) < \sqrt{2}\N<  \frac{z}b+c} \td z\\
&\quad +  \int_{\frac{2abc}{b-a}}^0 e^z \pk{  \frac{z}{a} -c < \sqrt{2}\N< \min\left(-\sqrt{c^2+\frac{2cz}{x}}, \frac{z}b+c\right)} \td z.
\end{align*}
Consequently, we give the explicit calculations of $H_3$ by the related claims of $H_{31}\sim H_{34}$.
}
\QED

\bibliographystyle{ieeetr}
\bibliography{reff}

\end{document}